\documentclass[12pt]{amsart}
\RequirePackage{newlfont}
\usepackage{amscd}

\usepackage{amsthm} 

\usepackage{amsmath}
\usepackage{amsfonts}
\usepackage{amssymb}
\usepackage{enumerate}
\usepackage{graphicx}
\usepackage{latexsym}
\usepackage{color}
\usepackage{bbm}

\usepackage{relsize} 

\newcommand{\bigint}{\mathlarger{\int}}

\newtheorem{theorem}{Theorem}[section]

\newtheorem{proposition}[theorem]{Proposition}
\newtheorem{lemma}[theorem]{Lemma}
\newtheorem{corollary}[theorem]{Corollary}

\newtheorem{definition}[theorem]{Definition}

\newtheorem{example}[theorem]{Example}
\newtheorem{remark}[theorem]{Remark}

\newtheorem{question}{Question}

\numberwithin{equation}{section}

\begin{document}

	\title[Modular topologies]{Modular Topologies on Vector Spaces}
	
	\author[M. A. Khamsi, J. Lang, O. M\'{e}ndez ]{Mohamed  A. Khamsi, Jan Lang,  Osvaldo M\'{e}ndez, }

	\address{Mohamed A. Khamsi\\ Department of Applied Mathematics and Sciences, Khalifa University, Abu Dhabi, UAE}
	\email{mohamed.khamsi@ku.ac.ae}
	\address{Jan Lang\\ Department of Mathematics, 100 Math Tower, 231 West 18th Ave., Columbus,
		OH 43210-1174, USA }
	\email{lang.162@osu.edu}
	\address{Osvaldo M\'{e}ndez\\Department of Mathematical Sciences, The University of Texas at El Paso, El Paso, TX 79968, USA}
	\email{osmendez@utep.edu}

	\subjclass[2020]{Primary 46A16, Secondary 46B20, 46E30, 46A20}
	\keywords{\(\Delta_2\)-property, Luxemburg norm, modular topology, modular vector space, modular convergence, variable exponent spaces.}
	
	\begin{abstract}

		This paper addresses the topological structures induced on vector spaces by convex modulars that do not satisfy the $\Delta_2$ condition, with  particular focus on their applications to variable exponent spaces such as \( \ell^{(p_n)} \) and \( L^{p(\cdot)} \). The motivation behind this investigation  is its applicability to the study of boundary value problems involving the variable exponent $p(x)$-Laplacian when $p(x)$ is unbounded, a line of research  recently opened by the authors. Fundamental topological properties are analyzed, including separation axioms, countability axioms, and the relationship between modular convergence and classical topological concepts such as continuity.  Attention is given to the relation between modular and norm topologies.
		Special emphasis is placed on the openness of modular balls, the impact of the \(\Delta_2\)-condition, and duality with respect to modular topologies. 
		
	\end{abstract}
	\maketitle

	\thispagestyle{empty}
	\section{Introduction and motivation}
	The modular structure was tacitly used by Riesz as early as 1910 in his seminal work \cite{riesz} and later formalized by Orlicz in \cite{orlicz1931}; see also \cite{birnbaum_orlicz_1931}. In 1950, Nakano \cite{nakano, nakano3} formalized and systematically studied the idea of a modular vector space. By then, the Banach space structure had already been established as a powerful tool for dealing with partial differential equations. At least in the $L^p$-theory for constant $p$, the modular and Banach space structures are topologically identical. It was only natural that PDE's specialists turned their attention to norms rather than to modulars.\\
	Later on, the need for a suitable mathematical framework for the description of the hydrodynamics of electro-rheological and magneto-rheological-fluids \cite{RaRu} naturally led to the utilization of variable-exponent Lebesgue spaces $L^{p(\cdot)}$, which was first systematically studied in \cite{KR}, see also \cite{ve-book, ML, sharapudinov, waterman}. It is well known that as long as the variable exponent $p(\cdot)$ remains bounded, modular convergence and norm convergence on $L^{p(\cdot)}$ are equivalent. For this reason, once again, the main attention was traditionally placed on the Banach space structure of the Lebesgue and the Sobolev spaces with variable exponent. This is, in some sense, counter-intuitive, since many of the differential equations of hydrodynamics are modular in nature. For example, for fixed $\theta\in W^{1,p(\cdot)}(\Omega)$ the problem of minimization of the Dirichlet integral on a domain $\Omega\subset {\mathbb R}^n$, over given subspace $V \subset W^{1,p(\cdot)}(\Omega):$
	\begin{equation}\label{minimization}
		W^{1,p(\cdot)}(\Omega)\supset V\ni u\rightarrow \int\limits_{\Omega}|\nabla (\theta -u)(x)|^{p(x)}dx
	\end{equation}
	is plainly a modular problem and certainly does not seem to be fully equivalent to a norm-minimization problem.\\
	The situation changes radically when the exponent $p(x)$ is unbounded. In this case, modular convergence on $L^{p(\cdot)}$ is weaker than norm convergence, and it engenders a topology, referred to as the modular topology, which is weaker than the norm-topology engendered by the usual Luxemburg norm. In a separate work \cite{AOJA}, the authors show that Banach space techniques are somewhat unsuitable to handle the problem (\ref{minimization}) when the function $p(x)$ is unbounded, whereas the problem is naturally posed and solved in terms of the underlying modular structure of $W^{1p(\cdot)}$.\\
	It is with these non-linear applications in mind that we present a detailed study of the modular topology associated with a convex modular and place special emphasis on the case of variable exponent Lebesgue spaces with unbounded variable exponent in order to get techniques required for a more detailed study of corresponding non-linear problems involving nonstandard growth.\\
	More precisely, the topological structure underlying a convex modular on a vector space will be thoroughly examined. First, the basic properties, such as separation, separability, first and second countability, and stability, with respect to the linear structure will be studied. It will be shown that modular balls are generally not open; in the case of the variable exponent Lebesgue spaces, the openness of modular balls is equivalent to the boundedness of the exponent, which in turn is equivalent to the $\Delta_2$-condition \cite{chen, kaminska, KK, kozlowski_book, rao}.\\
	The paper is structured as follow: In the next section the well known concept of modular vector space is introduced and general properties are recalled.\\
	In Section \ref{modulartopologies} the modular topology induced by a convex modular is studied, establishing its fundamental properties and its relationship with modular convergence. Key results in this section include the characterization of modularly open sets, the uniqueness of limits, and the interplay between modular and norm topologies; it is shown in particular that the validity of the \(\Delta_2\)-condition is equivalent to the modular topology being compatible with the linear structure of the underlying vector space.  
	Section \ref{lpn} provides a detailed investigation of the modular topology on the variable exponent spaces \(\ell^{(p_n)}\) and \(L^{p(\cdot)}\), emphasizing the structural differences it presents with normed spaces. A central result establishes that the openness of modular balls—and thus the equivalence of modular and norm topologies—is intrinsically linked to the boundedness of the exponent and the validity of the \(\Delta_2\)-condition. In the last section we study duality in modular vector spaces and characterize modularly continuous linear functionals. A key finding of this section is that in the absence of the $\Delta_2$-condition, modular duality  deviates from classical Banach space duality.

	\section{Modular vector spaces}\label{modularvectorspaces}
	Modular vector spaces provide a rich structure that allows for the exploration of aspects of functional analysis that are out of the reach of the theory of topological vector spaces. Unlike normed spaces, which are typically restricted by the linear properties of norms, the topology of modular spaces is generated by a convex modular function, offering more flexibility in addressing many peculiarities of non-linear problems. This flexibility makes them a powerful tool in areas such as optimization, duality theory, and approximation in areas involving non-linear phenomena.\\
	In this section, the definition and basic properties of modular vector spaces are recalled, starting with the notion of a convex modular. The goal is to build a foundation for understanding their topological properties, which will play a critical role in the rest of this work.
	
	\begin{definition}\label{modular-vs}
		A modular on a real vector space $X$ is a functional $\rho:X\rightarrow [0,\infty]$ such that
		\begin{enumerate}\label{modularproperties}
			\item[(1)] $\rho(u) = 0$ if and only if $u = 0$;
			\item[(2)] $\rho(\alpha u) = \rho(u)$, if $|\alpha| = 1$;
			\item[(3)] $\rho(\alpha u + (1-\alpha) v )\leq \alpha \rho(u) + (1-\alpha)\rho(v)$, for any $\alpha \in [0,1]$, and any $u, v \in X$.
		\end{enumerate}
	\end{definition}
	
	In what follows, we assume that $\rho$ is left-continuous, meaning that for any $u \in X$, $\lim\limits_{\lambda \rightarrow \lambda_0^-}\ \rho(\lambda u) = \rho(\lambda_0 u)$, for each $\lambda_0>0.$\\
	
	Let $X$ be a real vector space, and let $\rho:X \to [0,\infty]$ be a convex modular. Then,
	$$X_{\rho}=\Big\{v\in X: \rho(\lambda v)<\infty \, \text{for some} \, \lambda>0\Big\}$$
	is a vector subspace of $X$. It is straightforward to show that
	$$X_{\rho}=\left\{v\in X: \lim\limits_{\lambda \to 0} \rho(\lambda v) = 0\right\}.$$ 
	
	\noindent Two notable examples of modular vector spaces were introduced in the pioneering works of Orlicz \cite{orlicz1931} and Nakano \cite{nakano, nakano3}. Specifically:
	
	\begin{definition}\label{deflp}
		Fix a sequence $\mathbf{p}: 
		= (p_n) \subset [1,\infty)$. On the vector space $X = \mathbb{R}^{\mathbb{N}}$, define the functional $\rho_{\mathbf{p}}: X \to [0, \infty]$ by
		$$\rho_{\mathbf{p}}((a_j)) := \textstyle \sum\limits_{j = 1}^{\infty} |a_j|^{p_j}.$$
		Then $\rho_{\mathbf{p}}$ is a convex and left-continuous modular. Moreover, the associated modular vector space $X_{\rho_{\mathbf{p}}}$, denoted $\ell^{(p_n)}$, is 
		$$\ell^{(p_n)} := \Big\{(a_j) \in  \mathbb{R}^{\mathbb{N}} : \textstyle \sum\limits_{j = 1}^{\infty} |\lambda a_j|^{p_j} < \infty \, \text{for some} \, \lambda > 0 \Big\}.$$
		
	\end{definition}
	The following definition captures the continuous counterpart of the $\ell^{(p_n)}$ spaces. 
	\begin{definition}\label{Lp}
		\noindent Let $\Omega\subseteq {\mathbb R}^n$ be an open set and $p:\Omega\rightarrow [1,\infty)$ be a Borel-measurable function. On the set of extended-real valued Borel-measurable functions on $\Omega$, ${\mathcal M}(\Omega)$, the functional
		$$\rho_{p(\cdot)}(u):=\int\limits_{\Omega}|u(x)|^{p(x)}dx,$$ or shortly $\rho_p$,
		is a left-continuous convex modular; the associated modular space is
		$$L^{p(\cdot)}(\Omega):=\Bigg\{u\in {\mathcal M}(\Omega): \int\limits_{\Omega}|\lambda u(x)|^{p(x)}<\infty \,\,\text{for some} \lambda>0 \Bigg\}.$$
	\end{definition}
	\noindent In this work the modular vector spaces \(\ell^{(p_n)}\) and \(L^{p(\cdot)}(\Omega)\) will be referred to as variable exponent spaces.
	
	\begin{remark}{\normalfont  
			Consider a convex modular $\rho$ on a vector space $X$. Let 
			$$D:=\{x\in X_\rho:\ \rho(x)\leq 1\}.$$ 
			Then it holds $X_{\rho}=\langle D\rangle$, meaning $X_\rho$ is the linear span of $D$. Indeed, if $x \in X_{\rho}$, either $\rho(\lambda x) \leq 1$ or $\infty >\rho(\lambda x) > 1$. In the latter case, it holds that
			$$\rho\left(\frac{\lambda x}{\rho(\lambda x)}\right) \leq 1.$$
			In both cases, $x$ belongs to $\langle D\rangle$. Conversely, if $x\in \langle D\rangle$, let $x = \textstyle \sum\limits_{j=1}^K\alpha_jx_j$ with $x_j \in D$ for $j = 1, \dots, K$. Since $D$ is balanced, considering $\alpha_j > 0$ for all $j$ suffices. Then it follows from the convexity of $\rho$ that
			$$\rho\Big(\big(\textstyle \sum\limits_{j=1}^K \alpha_j\big)^{-1} x\Big) \leq 1,$$
			proving the claim.}
	\end{remark}
	
	The concept of modular balls is of key importance in what follows.
	
	\begin{definition}
		For $a \in X$ and $\varepsilon > 0$, let  
		$$B_{\rho,\varepsilon}(a) := \big\{y \in X : \rho(y - a) < \varepsilon\big\}.$$ 
		In the sequel $B_{\rho,\varepsilon}(a)$ will be referred to as the modular ball of radius $\varepsilon$, centered at $a$.
	\end{definition}
	
	\begin{proposition}\label{X_rho-stability}
		If $a \in X_{\rho}$, the convexity of $\rho$ implies that $B_{\rho,r}(a) \subseteq X_{\rho}$, for any $r > 0$.
	\end{proposition}
	\begin{proof}
		By assumption, there exists $\lambda > 0$ such that $\rho(\lambda a) < \infty$. Without loss of generality, assume $0 < \lambda < 1$. Then, for any $y \in B_{\rho,r}(a)$,
		\begin{align*}\rho\left(\frac{\lambda}{2} y\right) = \rho\left(\frac{\lambda}{2} (y - a) + \frac{\lambda}{2} a\right) &\leq \frac{1}{2} \Big( \rho(\lambda (y - a)) + \rho(\lambda a) \Big) \\
			&< \frac{1}{2} \big(\lambda\ r + \rho(\lambda a)\big) < \infty.
		\end{align*}
	\end{proof}
	
	\subsection{The Luxemburg norm}\leavevmode
	
	Since each modular ball $B_{r}(0):=\{x\in X_\rho: \rho(x) < r\}$ is convex, balanced, and absorbent set on $X_{\rho}$. The Minkowsky functional
	$$\mu_{B_r}(z) := \inf \Big\{\lambda > 0 : \rho(\lambda^{-1}z) \leq r \Big\}$$
	defines a norm on $X_{\rho}$. It can be easily shown that on account of the convexity of $\rho$, for any two modular balls $B_{r_1}(0)$ and $B_{r_2}(0)$ with $r_1 \leq r_2$, their corresponding Minkowsky functionals satisfy
	$$\mu_{B_{r_2}}(x) \leq \mu_{B_{r_1}}(x) \leq \frac{r_2}{r_1} \mu_{B_{r_2}}(x),$$
	for all $x \in X_{\rho}$. Due to this equivalence of norms, it is customary to focus on the Minkowsky functional of the unit modular ball, denoted as $\mu_{B_1}$. We also adopt the notation
	$$\mu_{B_1}(x) = \|x\|_{\rho} := \inf \Big\{\lambda > 0 : \rho(\lambda^{-1}x) \leq 1 \Big\}.$$
	The norm $\|\cdot\|_{\rho}$ is referred to as the Luxemburg norm on the modular space $X_{\rho}$ \cite{Luxemburg}.
	
	\begin{proposition}\label{standard}
		Consider the modular vector space $X_{\rho}$.
		\begin{enumerate}
			\item [$(i)$] Since $\rho$ is left-continuous, for any $x \neq 0$, then $\rho\left(\|x\|^{-1}_\rho x \right) \leq 1$.
			\item[$(ii)$] $\rho(x) \leq 1 \iff \|x\|_{\rho} \leq 1$.
			\item [$(iii)$] If $\|x\|_{\rho} \leq 1$, then $\rho(x) \leq \|x\|_{\rho}$.
			\item[$(iv)$] If, in addition, $\rho$ is right-continuous, then $\rho(x) < 1 \iff \|x\|_{\rho} < 1$.
		\end{enumerate}
	\end{proposition}
	\begin{proof}
		The proof follows directly from the definitions and properties of the modular $\rho$.
	\end{proof}
	
	\begin{remark}{\normalfont  
			In general, the condition $\rho(x_0) < 1$ does not imply $\|x_0\| < 1$. For instance, take the domain $\Omega = \left(0, \frac{1}{2}\right)$ and the function $p(x) = x^{-1} $. Consider the space $L^{p(\cdot)}((0,\frac{1}{2}))$ as in Definition \ref{Lp}, that is, the modular is given by  $\rho(u) = \bigint_0^{1/2} |u(x)|^{p(x)} dx$ .
			Let $\theta > 1$. Observe that 
			$$\lim_{x\rightarrow 0^+} \big(\theta^{x^{-1}} - x^{-1}\big) = \infty.$$
			Thus, there exists a point $x_0$ in the interval $0 < x_0 < \frac{1}{2}$ such that for $0 < x < x_0$, we have $\theta^{x^{-1}} - x^{-1} > \frac{1}{2}.$
			It follows that
			$$\int_0^{\frac{1}{2}} \theta^{x^{-1}} dx \geq \int_0^{x_0} \left( \frac{1}{2} + \frac{1}{x} \right) dx = \infty.$$
			Consequently, we have $\rho(1) = \frac{1}{2} < 1$. However, for any $\lambda$ such that $0 < \lambda < 1$, $\rho\left( \frac{1}{\lambda} \right) = \infty$, which implies that $\|1\|_{\rho} = 1$. This example illustrates that modular balls might not be norm-open.\\
			Moreover, this example demonstrates that equality may not always hold in $(i)$ if the modular fails to be right-continuous.}
	\end{remark}
	
	\section{Modular topologies}\label{modulartopologies}
	\ { In the sequel, a net $x:D\rightarrow X$ from a directed set $D$ into $X$ is said to $\rho$-converge to $x\in X$ iff for any $\delta>0$ there exists $\alpha_0\in D$ such that $\rho(x_{\alpha}-x)<\delta$ for $D\ni  \alpha \geq \alpha_0$. In particular, a sequence $(x_j)\subset X$ such that $\rho(x_j - x) \rightarrow 0$ as $j \rightarrow \infty$ $\rho$-converges to $x$. The notation $(x_{\alpha}) \overset{\rho}{\rightarrow} x$ will be used to denote $\rho$-convergence.} { Note that any net will $\rho$-converge to at most one point, meaning that the $\rho$-limit, if it exists, is unique. Indeed, let {  $(x_{\alpha}) \subset X$} $\rho$-converge to both $x$ and $y$. 
		Then, since
		$$\rho\left(\frac{x - y}{2}\right) \leq \frac{1}{2} \rho(x - x_{\alpha}) + \frac{1}{2} \rho(x_{\alpha} - y),$$
		letting $\alpha \to \infty$ gives $\rho\left(\frac{x - y}{2}\right) = 0$, which implies that $x = y$.} From Proposition \ref{X_rho-stability}, it follows that if {  $(x_{\alpha})$ is a net in $X$ and $(x_{\alpha}) \subset X_\rho$ $\rho$-converges to $x$, then $x \in X_\rho$.}\\
	
	The following result is straightforward:
	
	\begin{proposition}\label{rho-topology}
		Consider { $A\subset X$}. The following are equivalent:
		\begin{enumerate}
			\item For any sequence $(x_n)\subseteq A$ which $\rho$-converges to $x$, it holds $x \in A$;
			\item For any $x \notin A$, there exists $\varepsilon > 0$ such that $B_{\rho, \varepsilon}(x) \cap A = \emptyset$, i.e., $B_{\rho, \varepsilon}(x) \subset { X} \setminus A$.
		\end{enumerate}
	\end{proposition}
	
	The modular $\rho$ defines a natural topology {  on $X$}, denoted by $\tau_\rho$, as follows:
	
	\begin{definition}
		A subset $A$ of { $X$} is said to be $\tau_\rho$-open if for any $x \in A$, there exists $\varepsilon > 0$ such that $B_{\rho, \varepsilon}(x) \subset A$.
	\end{definition}
	
	From Proposition \ref{rho-topology}, it is clear that { $C \subset X$} is $\tau_\rho$-closed if and only if, for any sequence $(x_n)\subseteq C$ which $\rho$-converges to $x$, it holds $x \in C$.
	
	\begin{corollary}
		It follows from Proposition \ref{X_rho-stability} that the unique $\tau_{\rho}$-limit $x$ of a $\tau_{\rho}$-convergent sequence $(x_{n})\subset X_{\rho}$, must be in $X_{\rho}$. Therefore, $X_{\rho}$ is a $\tau_{\rho}$-closed subspace of $X$.
	\end{corollary} 
	
	\begin{remark}{\normalfont
			It is obvious that for any { $x \in X$}, the complement ${ X}\setminus \{x\}$ is $\tau_\rho$-open. Hence, $\tau_\rho$ is a $T_1$ topology.
		}
	\end{remark}
	
	The next theorem will be of fundamental importance:
	
	\begin{theorem}\label{convergencecharacterization}
		For any sequence $(x_n)$ in { $X$}, it holds that $(x_n) \overset{\rho}{\rightarrow} x$ if and only if $(x_n) \rightarrow x$ in $\tau_\rho$ (or shortly  $(x_n) \overset{\tau_\rho}{\rightarrow} x$).
	\end{theorem}
	\begin{proof}
		Assume first that $(x_n) \overset{\rho}{\rightarrow} x$. Let $\mathcal{O}$ be a $\tau_\rho$-open set that contains $x$. There exists $\varepsilon > 0$ such that $B_{\rho, \varepsilon}(x) \subset \mathcal{O}$. Since $(x_n) \overset{\rho}{\rightarrow} x$, there exists $n_0 \geq 1$ such that $x_n \in B_{\rho, \varepsilon}(x)$ for all $n \geq n_0$. Hence, $x_n \in \mathcal{O}$ for all $n \geq n_0$, that is, $(x_n) \overset{\tau_\rho}{\rightarrow} x$. Now assume that $(x_n) \overset{\tau_\rho}{\rightarrow} x$. Suppose that $(x_n)$ does not $\rho$-converge to $x$. Then, there exist $\varepsilon_0 > 0$ and a subsequence $(x_{n_k})$ such that $\rho(x - x_{n_k}) \geq \varepsilon_0$ for all $k \geq 1$. The subsequence $(x_{n_k})$ has a $\rho$-convergent subsequence or it does not have a subsequence $\rho$-converges. In the first case, let $(y_j)$ be a subsequence of $(x_{n_k})$ such that $y_j \overset{\rho}{\rightarrow} y \in X_\rho$. Clearly, $x \neq y$. Then, the set
		$$C = \{y_j : j \in \mathbb{N}\} \cup \{y\}$$
		is $\tau_\rho$-closed and $x \notin C$. Since $(y_j) \overset{\tau_\rho}{\rightarrow} x$, we have $x \in C$, which is a contradiction. On the other hand, if no subsequence of $(x_{n_k})$ $\rho$-converges to any point, the set
		$$A = \{x_{n_k} : k \in \mathbb{N}\}$$
		is $\tau_\rho$-closed in $X_\rho$ and does not contain $x$. Its complement, $X_\rho \setminus A$, is thus $\tau_\rho$-open and contains $x$. Hence, $X_\rho \setminus A$ must contain $x_{n_k}$ for sufficiently large $n_k$, leading to a contradiction.
	\end{proof}
	
	Theorem \ref{convergencecharacterization}, along with the uniqueness of the \(\rho\)-limit, gives the following result:
	
	\begin{corollary}
		Any sequence in { $X$} can $\tau_{\rho}$-converge to at most one limit.
	\end{corollary}
	{ \begin{corollary}\label{nets}
			For any net $(x_{\alpha})\subset { X}$the following statements are equivalent:
			\begin{itemize}
				\item [$(i)$] $x_{\alpha}\overset{\tau_{\rho}}{\longrightarrow}x$\\
				\item[$(ii)$] $x_{\alpha}\overset{\rho}{\longrightarrow}x$.
			\end{itemize}
		\end{corollary}
		\begin{proof}
			The proof of the implication $(ii)\Rightarrow (i)$ is elementary and will be omitted.\\ 
			
			We show $(i) \Rightarrow (ii)$ via contradiction. Let $D$ be a directed set and $x:D\rightarrow X$ be a net. Assume $(i)$ and suppose there exists $\delta>0$ and a sequence $\alpha_1<\alpha_2<,...<\alpha_k....$ such that $\rho(x-x_{\alpha_k})\geq \delta.$ It is easy to verify that the sequence $(x_{\alpha_k})$ is a subnet of $(x_{\alpha})$ and therefore it must $\tau_{\rho}$-converge to $x$. Theorem \ref{convergencecharacterization} yields $x_{\alpha_k}\overset{\rho}{\rightarrow}x$, which is not possible. This contradiction proves the claim.
	\end{proof}}
	{ \begin{corollary}
			Any net in { $X$} can $\tau_{\rho}$-converge to at most one limit.
		\end{corollary}
		
		Note that $\tau_{\rho}$ is a topology on $X$ and then by  $\overline{\tau_{\rho}}$ we denote the subspace topology induced on $X_{\rho}$ by $\tau_{\rho}$. 
		\\On account of Proposition \ref{X_rho-stability} it can be quickly seen that  $A\subseteq X_{\rho}$ is $\overline{\tau_{\rho}}$-open iff for any $a\in A$ there exists $\delta>0$ such that  $\{y\in X_{\rho}:\rho(a-y)<\delta\}\subset A$.\\
		Thus:
		\begin{corollary}
			It holds: $\overline{\tau_{\rho}}=2^{X_{\rho}}\cap \tau_{\rho}$.
		\end{corollary}
		From now on, all statements refer to the topological space $(X_{\rho},\overline{\tau_{\rho}})$. By virtue of the preceding remark, any subset of $X_{\rho}$ is $\overline{\tau_{\rho}}$-open if and only if it is $\tau_{\rho}$-open. In the interest of notational simplicity, we will slightly abuse the notation and use  $\tau_{\rho}$ instead of $\overline{\tau_{\rho}}$ to denote the subspace topology on $X_{\rho}$.
	}
	\begin{definition}
		Let $\overline{A}^\rho$ denote the {  $\tau_{\rho}$- closure of a set $A\subseteq X_{\rho}$} i.e.
		\[
		\overline{A}^\rho := \bigcap\limits_{A \subseteq W\ \text{$\tau_\rho$-closed}} W.
		\]
	\end{definition}

	The following properties will be used frequently:

	\begin{proposition}\label{basic-properties}
		The following properties hold:
		\begin{enumerate}
			\item If $\mathcal{O}$ is a $\tau_\rho$-open subset of $X_\rho$, then $\mathcal{O} + x = \{u + x; u \in \mathcal{O}\}$ is also $\tau_\rho$-open, for any $x \in X_\rho$. Hence, $\mathcal{O}_1 + \mathcal{O}_2$ is $\tau_\rho$-open provided either $\mathcal{O}_1$ or $\mathcal{O}_2$ is $\tau_\rho$-open.
			\item If $\mathcal{O}$ is $\tau_\rho$-open and $\alpha \geq 1$, then $ \alpha  \mathcal{O}$ is also $\tau_\rho$-open.
			\item For any $x \in \overline{A}^\rho$ and any $\tau_\rho$-open subset $\mathcal{O}$ such that $x \in \mathcal{O}$, we have $\mathcal{O} \cap A \neq \emptyset$.
			\item If $A$ is convex, then $\overline{A}^\rho$ is convex.
		\end{enumerate}
	\end{proposition}
	\begin{proof}\leavevmode
		\begin{enumerate}
			\item Let $y \in \mathcal{O} + x$, then $y - x \in \mathcal{O}$. Since $\mathcal{O}$ is $\tau_\rho$-open, there exists $\varepsilon > 0$ such that $B_{\rho,\varepsilon}(y - x) \subset \mathcal{O}$. Clearly, $B_{\rho,\varepsilon}(y) \subset \mathcal{O} + x$. As for $\mathcal{O}_1 + \mathcal{O}_2$, note that 
			\[\mathcal{O}_1 + \mathcal{O}_2 = \bigcup_{y \in \mathcal{O}_2} \mathcal{O}_1 + y = \bigcup_{x \in \mathcal{O}_1} \mathcal{O}_2 + x,\]
			yielding the desired conclusion.
			\item Let $x \in  \alpha  \mathcal{O}$. Then $ \beta  x \in \mathcal{O}$, where $\beta = 1/\alpha \in (0,1]$. Since $\mathcal{O}$ is $\tau_\rho$-open, there exists $\varepsilon > 0$ such that $B_{\rho,\varepsilon}( \beta  x) \subset \mathcal{O}$. For any $y \in B_{\rho,\alpha \varepsilon}(x)$, we have
			\[\rho\left( \beta  x -  \beta  y\right) \leq  \beta  \rho(x - y) <  \beta   \alpha  \varepsilon = \varepsilon,\]
			i.e., $ \beta  y  \in B_{\rho,\varepsilon}( \beta  x) \subset \mathcal{O}$. Hence $y \in  \alpha  \mathcal{O}$, which forces $B_{\rho, \alpha \varepsilon}(x) \subset \alpha \ \mathcal{O}$. Therefore, $\alpha \ \mathcal{O}$ is $\tau_\rho$-open. 
			\item Suppose $\mathcal{O} \cap A = \emptyset$. Then $A \subset \mathcal{O}^c = X_\rho \setminus \mathcal{O}$. Since $\mathcal{O}$ is $\tau_\rho$-open, it follows that $\mathcal{O}^c$ is $\tau_\rho$-closed. Hence, $\overline{A}^\rho \subset \mathcal{O}^c$, i.e., $\mathcal{O} \cap \overline{A}^\rho = \emptyset$. This contradicts the assumption that $x \in \mathcal{O} \cap \overline{A}^\rho$.
			\item Suppose $A$ is convex. Let $x, y \in \overline{A}^\rho$ and $\alpha \in (0, 1)$.  Assume $\alpha x + (1 - \alpha) y$ is not in $\overline{A}^\rho$.  Since $\mathcal{O} = X_\rho \setminus \overline{A}^\rho$ is $\tau_\rho$-open and $\alpha x + (1 - \alpha)y \in \mathcal{O}$, we can deduce that
			\[x \in \mathcal{O}_1 =  \beta  \mathcal{O} - (\beta -1)\ y,\]
			where $\beta = 1/\alpha$.  Since $\mathcal{O}_1$ is $\tau_\rho$-open, by the previously proven properties, $A \cap \mathcal{O}_1 \neq \emptyset$. Let $a \in A \cap \mathcal{O}_1$. Then there exists $x_0 \in \mathcal{O}$ such that
			\[a =  \beta  x_0 - (\beta -1)\ y,\]
			implying
			\[y = \frac{1}{1 - \alpha} x_0 - \frac{\alpha}{1 - \alpha} a \in \mathcal{O}_2 = \frac{1}{1 - \alpha} \mathcal{O} - \frac{\alpha}{1 - \alpha} a.\]
			Thus, $A \cap \mathcal{O}_2 \neq \emptyset$. Let $b \in A \cap \mathcal{O}_2$. Since $b \in \mathcal{O}_2$, there exists $y_0 \in \mathcal{O}$ such that
			\[b = \frac{1}{1 - \alpha} y_0 - \frac{\alpha}{1 - \alpha} a,\]
			implying $y_0 = \alpha a + (1 - \alpha)b \in A$, by the convexity of $A$. Therefore, $y_0 \in A \cap \mathcal{O}$ which contradicts the assumption that $\mathcal{O} \cap A = \emptyset$.
		\end{enumerate}
	\end{proof}
	
	The preceding result has the following fundamental consequence:
	
	\begin{proposition}\label{subspace}
		Let $A$ be a vector subspace of $X_\rho$. Then $\overline{A}^\rho$ is a $\tau_\rho$-closed vector subspace of $X_\rho$.
	\end{proposition}
	\begin{proof}
		It will be shown first that if $x, y \in \overline{A}^\rho$, then $x + y \in \overline{A}^\rho$. Suppose not, i.e., $x + y \in \mathcal{O} = X_\rho \setminus \overline{A}^\rho$. Then $x \in (\mathcal{O} - y)$, which is $\tau_\rho$-open. This forces $A \cap (\mathcal{O} - y) \neq \emptyset$. Let $a \in A \cap (\mathcal{O} - y)$. There exists $x_0 \in \mathcal{O}$ such that $a = x_0 - y$, implying $y = x_0 - a \in (\mathcal{O} - a)$. Again, since $\mathcal{O} - a$ is open, it follows that $A \cap (\mathcal{O} - a) \neq \emptyset$. Let $b \in A \cap (\mathcal{O} - a)$. Then there exists $y_0 \in \mathcal{O}$ such that $b = y_0 - a$, implying $y_0 = b + a$. Since $A$ is a subspace, we have $b + a \in A$. Therefore, the assumption implies that $y_0 \in A \cap \mathcal{O}$, contradicting the fact that $A \cap \mathcal{O} = \emptyset$. \\
		Next, observe that if $x \in \overline{A}^\rho$ and $\alpha \in \mathbb{R}$, then $ \alpha  x \in \overline{A}^\rho$. Without loss of generality, assume $\alpha \neq 0$.  Take first $\alpha > 0$. The first part of the proof shows that $k\ x \in \overline{A}^\rho$ for any $k \in \mathbb{N}$. Thus, it may be assumed that $\alpha$ is not an integer. In this case, there exists $k \in \mathbb{N}$ such that $k < \alpha < k + 1$. This implies the existence of $\theta \in (0,1)$ such that $\alpha = \theta k + (1 - \theta)(k + 1)$. Hence,
		\[ \alpha  x = \theta\ k\ x + (1 - \theta)(k + 1)\ x.\]
		Hence, $\alpha \ x \in \overline{A}^\rho$ since $\overline{A}^\rho$ is convex by Proposition \ref{basic-properties}.\\
		Finally, we show that if $x \in \overline{A}^\rho$, then $-x \in \overline{A}^\rho$. This follows from similar reasoning and the fact that if $\mathcal{O}$ is $\tau_\rho$-open, then $-\mathcal{O}$ is also $\tau_\rho$-open. To see this, let $\mathcal{O}$ be a $\tau_\rho$-open subset of $X_\rho$. Let $y \in -\mathcal{O}$. Then $-y \in \mathcal{O}$, implying the existence of $\varepsilon > 0$ such that $B_{\rho,\varepsilon}(-y) \subset \mathcal{O}$. Using the properties of the modular, it is easily seen that $z \in B_{\rho,\varepsilon}(-y)$ if and only if $-z \in B_{\rho,\varepsilon}(y)$. Therefore, $B_{\rho,\varepsilon}(y) \subset -\mathcal{O}$, completing the proof that $-\mathcal{O}$ is $\tau_\rho$-open.
	\end{proof}
	
	\subsection{The modular topologies $\tau_\lambda$}\leavevmode
	
	For each $\lambda > 0$, define \(\rho_{\lambda}: X \to [0, \infty]\) by \(\rho_{\lambda}(x) = \rho(\lambda x)\). It is straightforward to see that \(\rho_{\lambda}\) is a convex modular on \(X\). Moreover, for any \(\lambda > 0\), we have
	\[
	X_{\rho_{\lambda}} := \{x \in X: \rho_{\lambda}(\alpha x) < \infty \, \text{for some}\, \alpha > 0\} = X_{\rho},
	\]
	i.e., all the modulars \(\rho_{\lambda}\) define the same modular vector space \(X_{\rho}\). Let \(\tau_{\lambda}\) denote the \(\rho_{\lambda}\)-modular topology on \(X_{\rho}\). If \(0 < \lambda_1 < \lambda_2\), then any \(\tau_{\lambda_1}\)-closed set in \(X_{\rho_{\lambda_2}}\) is also \(\tau_{\lambda_2}\)-closed. Hence, \(\tau_{\lambda_1}\) is weaker than \(\tau_{\lambda_2}\). Simply put, the family of topologies \((\tau_{\lambda})\) increases with \(\lambda\).\\
	
	\begin{definition}
		Let \(\rho\) be a convex modular on a vector space \(X\).
		\begin{enumerate}
			\item The final topology of the family \((\tau_{\lambda})\) is given by
			\[
			\tau^f := \bigcap_{\lambda > 0} \tau_{\lambda}.
			\]
			\item The initial topology, denoted \(\tau^i\), is the weakest topology that is stronger than all \(\tau_{\lambda}\) for \(\lambda > 0\).
		\end{enumerate}
	\end{definition}
	
	The open subsets for both topologies \(\tau^f\) and \(\tau^i\) are characterized in the following proposition:
	
	\begin{proposition}\label{tau^i-tau^f-characterization}
		The following hold:
		\begin{enumerate}
			\item The open sets in \(\tau^f\) are those subsets \(Y \subseteq X_{\rho}\) that are \(\tau_{\lambda}\)-open for every \(\lambda > 0\).
			\item The open sets in \(\tau^i\) consist of arbitrary unions of finite intersections of the form \(A_{\lambda_1} \cap A_{\lambda_2} \dots \cap A_{\lambda_N}\), where \(A_{\lambda_j} \in \tau_{\lambda_j}\). By virtue of the inclusion \(\tau_{\alpha} \subseteq \tau_{\beta} \) whenever \(\alpha \leq  \beta \), given any open set $A$ in \(\tau^i\) there exist $J\subseteq (0,\infty)$ such that
			\[
			A=\bigcup_{\lambda \in J \subseteq (0, \infty)} A_{\lambda},
			\]
			where each \(A_{\lambda}\) is open in \(\tau_{\lambda}\).
		\end{enumerate}
	\end{proposition}
	
	\begin{remark}\label{open-charac-tau^i}
		Let \(\mathcal{O}\) be an open subset in \(\tau^i\). According to Proposition \ref{tau^i-tau^f-characterization}, there exists a subset \(J \subset (0, \infty)\) and a family of subsets \((A_{\lambda})_{\lambda \in J}\) such that \(A_{\lambda}\) is open in \(\tau_{\lambda}\) for every \(\lambda \in J\), and
		\[
		\mathcal{O} = \bigcup_{\lambda \in J \subseteq (0, \infty)} A_{\lambda}.
		\]
		Without loss of generality, assume \(J\) is nonempty. Define \(\mathcal{O}_n := \bigcup_{\lambda \in J \cap (0, n]} A_{\lambda}\) for any \(n >0\). Then \(\mathcal{O}_n\) is open in \(\tau_n\) and \(\mathcal{O}_n \subset \mathcal{O}_{n+1}\) for all \(n >0 \). Moreover, it is easy to verify that
		\[
		\mathcal{O} = \bigcup_{n >0} \mathcal{O}_n.
		\]
	\end{remark}
	
	The next two propositions describe the initial and final topologies corresponding to the family \((\tau_{\lambda})\). Recall that for a family \(\mathcal{F}\) of topologies on a set \(X_{\rho}\), the initial topology of \(\mathcal{F}\) is the weakest topology on \(X_{\rho}\) that is stronger than each member of \(\mathcal{F}\), whereas the final topology is the strongest topology on \(X_{\rho}\) contained in each member of \(\mathcal{F}\).
	
	It is a routine exercise to verify that a set \(C \subseteq X_{\rho}\) is \(\tau^f\)-closed if and only if \(C\) is \(\tau_{\lambda}\)-closed for every \(\lambda > 0\), i.e., given any sequence \((x_j) \subseteq C\) such that, for some \(\lambda > 0\), \(\rho(\lambda (x_i - x)) \to 0\) as \(j \to \infty\), we have \(x \in C\). Since \(\tau^f\)-open subsets are also \(\tau_{\lambda}\)-open for every \(\lambda > 0\), we obtain the following fact:
	
	\begin{proposition}
		If $(x_j)\subset X_{\rho}$ and there exists $\lambda>0$ such that $(x_j)\overset{\tau_\lambda}{\rightarrow}x$, then $(x_j)\overset{\tau^f}{\rightarrow}x$.
	\end{proposition}
	
	Some algebraic properties of $\tau^i$ will be discussed next:
	
	\begin{proposition}\label{tau^i-algebra}
		The topology $\tau^{i}$ on $X_\rho$ is stable under addition and scalar multiplication. Specifically, if $A\in \tau^{i}$, $B\in \tau^{i}$, and $r\in {\mathbb R}$, with $r \neq 0$, then $A+B\in \tau^{i}$ and $rA\in \tau^{i}$.
	\end{proposition}
	\begin{proof}
		Observe first that if $A\in \tau^{i}$ and $B\in \tau^{i}$, then $A+B\in \tau^{i}$. For, according to Remark \ref{open-charac-tau^i}, $A=\bigcup_{n \geq 1}A_n$ and $B=\bigcup_{n \geq 1}B_n$, where $A_n$ and $B_n$ are in $\tau_n$ for all $n \geq 1$. Using Proposition \ref{basic-properties}, we know that $A_n + B_n \in \tau_n$ for all $n \geq 1$. It follows that $A + B = \bigcup_{n \geq 1}(A_n + B_n)$, which proves the desired result.
		
		On the other hand $rA$ is $\tau^{i}$-open for any $r \in \mathbb{R}$, $r \neq 0$, provided that $A$ is $\tau^{i}$-open. To see this, assume $r > 0$. Since
		\[
		\rho_{\frac{\lambda}{r}}(x-w)= \rho\left(\frac{\lambda}{r} (x - w)\right) = \rho_\lambda \left(\frac{x}{r} - \frac{w}{r}\right),
		\]
		it follows that $r B_{\rho_\lambda, \delta}\left(\frac{x}{r}\right) = B_{\rho_{\lambda/r}, \delta}(x)$ for modular balls, for any $\delta > 0$. Hence, if $A \in \tau_{\lambda}$ for some $\lambda > 0$, then $rA \in \tau_{\lambda/r}$. The properties of the family $(\tau_\lambda)_{\lambda > 0}$, guarantee that if $A \in \tau_n$ for some $n \geq 1$, then there exists $m \geq 1$ such that $rA \in \tau_m$. Since $\rho$-balls are symmetric, this result extends to any $r \neq 0$. Finally, using the characterization of $\tau^{i}$-open subsets, it is readily seen that $rA$ is in $\tau^i$ for any $r \neq 0$, which completes the proof of Proposition \ref{tau^i-algebra}.
	\end{proof}
	
	The next proposition establishes a comparison between the topology $\tau^i$ and the Luxemburg norm topology.
	
	\begin{proposition}\label{sequential-comparison}
		For any sequence $(x_j) \subset X_\rho$, the following are equivalent:
		\begin{enumerate}
			\item[(i)] $(x_j) \overset{\tau^i}{\rightarrow} x$;
			\item[(ii)] $(x_j) \overset{\tau_\lambda}{\rightarrow} x$, for all $\lambda > 0$.
		\end{enumerate}
	\end{proposition}
	\begin{proof}
		Without loss of generality, assume \(x = 0\). Suppose $(i)$ holds but $(ii)$ fails. Then there exists $\lambda_0 > 0$ such that $\rho_{\lambda_0}(x_j) \nrightarrow 0$ as $j \rightarrow \infty$. Hence, there exists $\varepsilon_0 > 0$ and a subsequence $(x_{n_k})$ such that $\rho_{\lambda_0}(x_{n_k}) \geq \varepsilon_0$. Denote the subsequence $(x_{n_k})$ by $(y_k)$. There are two possibilities for the sequence $(y_k)$: either it has a subsequence that $\rho_{\lambda_0}$-converges to a point $y \in X_\rho$ (which must be different from $0$), or there is no $\rho_{\lambda_0}$-convergent subsequence of $(y_k)$. 
		
		In the first case, the set $W = \{y_k, k \in {\mathbb N}\} \cup \{y\}$ is $\rho_{\lambda_0}$-closed and does not contain $0$. Since $\tau^i$ contains $\tau_{\lambda_0}$, $W$ is $\tau^i$-closed. However, since $(y_k) \overset{\tau^i}{\rightarrow} 0$, we must conclude that $0 \in W$, which is a contradiction. 
		
		In the second case, assume that $(y_k)$ has no $\rho_{\lambda_0}$-convergent subsequence. In this case, the set $G = \{y_k, k \in {\mathbb N}\}$ is $\rho_{\lambda_0}$-closed and does not contain $0$. Again, since $\tau_{\lambda_0}$-closed subsets are $\tau^i$-closed, we conclude that $G$ is $\tau^i$-closed, which implies that $0 \in G$ (since $(y_k) \overset{\tau^i}{\rightarrow} 0$), which is again a contradiction. Therefore, $(i)$ implies $(ii)$.
		
		Conversely, assume $\rho_{\lambda}(x_j) \rightarrow 0$ as $j \rightarrow \infty$ for all $\lambda > 0$. Let $V$ be an arbitrary $\tau^i$-neighborhood of $0$. Then there exists a $\tau^i$-open set $\mathcal{O}$ such that $0 \in \mathcal{O}$. According to Remark \ref{open-charac-tau^i}, we have ${\mathcal O} = \bigcup_{n \geq 1} A_n$, where each $A_n$ is $\tau_n$-open for any $n \geq 1$. Since $0 \in \mathcal{O}$, there exists some $n \geq 1$ such that $0 \in A_n$. Since $(x_j) \overset{\tau_n}{\rightarrow} 0$ and $A_n$ is a $\tau_n$-neighborhood of $0$, there exists $j_0 \geq 1$ such that for all $j \geq j_0$, we have $x_j \in A_n \subset \mathcal{O} \subset V$. Therefore, $(x_j) \overset{\tau^i}{\rightarrow} 0$, which completes the proof of Proposition \ref{sequential-comparison}.
	\end{proof}
	
	\noindent Proposition \ref{sequential-comparison} implies, in particular, that any \(\tau^i\)-convergent sequence has a unique limit within \(X_\rho\).
	
	Next, the Luxemburg norms associated to the modulars $\rho_\lambda$ for $\lambda > 0$ are considered.
	
	\begin{proposition}\label{equivalence}
		Consider the modular vector space $X_\rho$.
		\begin{itemize}
			\item[(i)] For any $\alpha > 0$ and $x \in X_\rho$, we have
			\[
			\begin{cases}
				\alpha \|x\|_\rho \leq \|x\|_{\rho_\alpha} \leq \|x\|_{\rho}, & \text{if } \alpha \leq 1, \\
				\|x\|_\rho \leq \|x\|_{\rho_\alpha} \leq \alpha \|x\|_{\rho}, & \text{if } \alpha > 1.
			\end{cases}
			\]
			\item[(ii)] The Luxemburg norms generated by $\rho_{\alpha}$ and $\rho_{\beta}$ for $\alpha, \beta > 0$ are equivalent (and hence they induce the same topology).
			\item[(iii)] Convergence with respect to \(\|\cdot\|_{\rho}\) in \(X_{\rho}\) is equivalent to $\rho_{\lambda}$-convergence for all $\lambda > 0$.
		\end{itemize}
	\end{proposition}
	\begin{proof}
		The proofs of $(i)$ and $(ii)$ are straightforward. As to the proof of $(iii)$, it will be next shown that $(x_n) \overset{\rho_\lambda}{\rightarrow} 0$ for all $\lambda > 0$ if and only if $(x_n) \overset{\|\cdot\|_\rho}{\rightarrow} 0$. 
		
		Assume $(x_n) \overset{\|\cdot\|_\rho}{\rightarrow} 0$. Fix $\lambda > 0$. For large $n$ it is clear that $\|\lambda x_n\| \leq 1$. Then it follows, on account of Proposition \ref{standard}, that $\rho(\lambda x_n) \leq \|\lambda x_n\| = \lambda \|x_n\|$, which implies $(x_n) \overset{\rho_\lambda}{\rightarrow} 0$. 
		
		Conversely, assume $(x_n) \overset{\rho_\lambda}{\rightarrow} 0$ for all $\lambda > 0$, and that $(x_n)$ does not converge to $0$ with respect to the norm $\|\cdot\|_\rho$. Without loss of generality, it may be assumed that there exists $\varepsilon_0 > 0$ such that $\|x_n\|_\rho > \varepsilon_0$ for all $n \in \mathbb{N}$. The definition of the Luxemburg norm yields $\rho(x_n / \varepsilon_0) > 1$ for all $n \in \mathbb{N}$. Hence, $(\rho_{1/\varepsilon_0}(x_n))$ will not converge to $0$, contradicting our assumption.
	\end{proof}
	
	\begin{theorem}
		Consider the collection $\tau_*$ consisting of those sets $A \subset X_{\rho}$ that satisfy the following property:\\
		For any $x \in A$, there exist $\lambda, \varepsilon > 0$ (both depending on $x$) such that the modular ball $B_{\rho_{\lambda}, \varepsilon}(x) \subseteq A$.\\
		Then the following hold:
		\begin{itemize}
			\item[(i)] \(\tau_*\) is a topology that contains all topologies \(\tau_{\lambda}\) for \(\lambda > 0\).
			\item[(ii)] If $(x_k)$ is any sequence in $X_{\rho}$, then $(x_k)\overset{\tau_{\ast}}{\rightarrow} x$ if and only if $(x_k)\overset{\rho_{\lambda}}{\rightarrow} x$ for all $\lambda > 0$.
		\end{itemize}
		Let $\tau_{\|\cdot\|_{\rho}}$ denote the topology induced on $X_{\rho}$ by the Luxemburg norm. Then $\tau_{\|\cdot\|_{\rho}} = \tau_*$, and $\tau_{\|\cdot\|_{\rho}}$ is the weakest first-countable topology that is stronger than each $\tau_{\lambda}$ for all $\lambda > 0$.
	\end{theorem}
	
	\begin{proof}
		For $(i)$, it is straightforward to show that $\tau_*$ is a topology. Indeed, if $A \in \tau_*$ and $B \in \tau_*$, there must exist $\lambda, \beta, \varepsilon_1, \varepsilon_2 \in (0,\infty)$ such that $B_{\rho_{\lambda}, \varepsilon_1}(x) \subseteq A$ and $B_{\rho_{\beta}, \varepsilon_2}(x) \subseteq B$. The properties of the modulars $(\rho_\alpha)$ imply that the modular ball $B_{\rho_{\beta}, \min\{\varepsilon_1, \varepsilon_2\}}(x)$ is contained in $A \cap B$. The verification of the fact that $\tau_*$ is closed under arbitrary unions is straightforward.
		
		It is evident by definition that for any $\lambda > 0$, $\tau_{\lambda} \subseteq \tau_*$.
		
		To prove $(ii)$, assume that $(x_k) \overset{\tau_{\ast}}{\rightarrow} x$. Suppose there exists $\lambda_0 > 0$ such that $\rho_{\lambda_0}(x - x_k) \not\rightarrow 0$ as $k \rightarrow \infty$. Then there exists $\delta > 0$ and a subsequence $(x_{k_j})$ such that $\rho_{\lambda_0}(x - x_{k_j}) \geq \delta$ for all $j \in \mathbb{N}$. Set $S = \{x_{k_j}; j \in \mathbb{N}\}$. Clearly, $x \notin S$.
		
		The subsequence $(x_{k_j})$ either contains a $\rho_{\lambda_0}$-convergent subsequence or no subsequence of $(x_{k_j})$ $\lambda_0$-converges. In the first case, select a subsequence, say $(y_i)$, such that $\rho_{\lambda_0}(y_i - y) \rightarrow 0$ as $i \rightarrow \infty$. Set $B = \{y_i; i \in \mathbb{N}\} \cup \{y\}$. Necessarily, $x \notin B$. The set $B$ is $\rho_{\lambda_0}$-closed in $X_{\rho}$ and therefore its complement $X_{\rho} \setminus B$ is $\rho_{\lambda_0}$-open, hence $\tau_{\ast}$-open, and it contains $x$. But $(y_i)$ is a subsequence of $(x_k)$, which $\tau_{\ast}$-converges to $x$. This is clearly a contradiction.
		
		Similarly, if no subsequence of $(x_{k_j})$ is $\rho_{\lambda_0}$-convergent, then $S$ is $\rho_{\lambda_0}$-closed, hence $\tau_{\ast}$-closed, and $S$ does not contain $x$. Again, a contradiction is reached by observing that in this case, $X_{\rho} \setminus S$ is a $\tau_{\ast}$-open set containing $x$.
		
		Conversely, assume $(x_k) \overset{\rho_{\lambda}}{\rightarrow} x$ for all $\lambda > 0$. Let $V$ be a neighborhood of $0$ in $\tau_{\ast}$. By definition, there exist $\delta > 0$, $\varepsilon > 0$ such that the modular ball $B_{\rho_{\delta}, \varepsilon}(0) \subset V$. Since $\rho_{\delta}(x - x_j) \rightarrow 0$ as $j \rightarrow \infty$ it is immediate that $x - x_j \in B_{\rho_{\delta}, \nu}(0)$ for large enough $j$. Thus, $x \overset{\tau_{\ast}}{\rightarrow} x$, as claimed.
		
		The inclusion $\tau_{\|\cdot\|_{\rho}} \subseteq \tau_{\ast}$ is tackled next. To this end, let $A \subseteq X_{\rho}$ be $\|\cdot\|_{\rho}$-open and $x \in A$. Then there exists $\varepsilon > 0$ such that $\{y \in X_\rho : \|y - x\|_{\rho} < \varepsilon\} \subset A$. By the definition of the Luxemburg norm, if $\rho_{\frac{2}{\varepsilon}}(x - y) < 1$, then $\|y - x\|_{\rho} < \varepsilon$. Hence, for any $x \in A$, the modular ball $B_{\rho_{\frac{\varepsilon}{2}}, 1}(x) \subset A$. It follows that $A$ is $\tau_{\ast}$-open.
		
		On the other hand, let $V$ be a $\tau_{\ast}$-open set. Take $x \in V$. By definition, there exist $\lambda_x > 0$ and $\varepsilon_x < 1$ such that 
		$$
		\rho(\lambda_x(y - x)) < \varepsilon_x \Rightarrow y \in V.
		$$
		If $\|\lambda_x(x - y)\|_{\rho} < \varepsilon_x < 1$, then one has, by virtue of Proposition \ref{standard} (iii), that
		\begin{align*}
			\rho(\lambda_x(y - x)) &= \rho\Big(\|\lambda_x(x - y)\|_{\rho}\|\lambda_x(x - y)\|^{-1}_{\rho}\lambda_x(y - x)\Big) \\
			&\leq \|\lambda_x(x - y)\|_{\rho} \rho\Big(\|\lambda_x(x - y)\|^{-1}_{\rho}\lambda_x(x - y)\Big) \\
			&< \varepsilon_x.
		\end{align*}
		It follows that the norm ball $\{y : \|x - y\|_{\rho} < \varepsilon_x \lambda_x^{-1}\}$ is contained in $V$, and thus $V$ is open in $\tau_{\|\cdot\|_{\rho}}$. Hence, $\tau_{\ast} \subseteq \tau_{\|\cdot\|_{\rho}}$, as claimed.
		
		Let $\tau$ be a first-countable topology on $X_{\rho}$ that is stronger than every $\tau_{\lambda}$ for $\lambda > 0$. By assumption, the inclusion $(X_{\rho}, \tau) \hookrightarrow (X_{\rho}, \tau_{\lambda})$ is continuous for every $\lambda > 0$. Thus, for any sequence $(x_j)$, the condition $x_j \overset{\tau}{\rightarrow} x$ implies $x_j \overset{\tau_{\lambda}}{\rightarrow} x$ for every $\lambda > 0$. But by Proposition \ref{equivalence}, this is equivalent to $x_j \overset{\|\cdot\|_{\rho}}{\rightarrow} x$. First countability implies that the inclusion
		\[
		(X_{\rho}, \tau) \hookrightarrow (X_{\rho}, \tau_{\|\cdot\|_{\rho}})
		\]
		is continuous, which means that $\tau_{\|\cdot\|_{\rho}}$ is weaker than $\tau$.
	\end{proof}
	
	The question arises whether the modular topology is compatible with the algebraic structure of \(X_{\rho}\). To clarify this, the following definition is necessary.
	
	\begin{definition}
		The modular $\rho$ is said to satisfy the $\Delta_2$-property if any sequence $(x_j) \subset X$ such that $(x_j)\overset{\rho}{\rightarrow}0$ also satisfies $(2x_j)\overset{\rho}{\rightarrow}0$.
	\end{definition}
	
	\begin{theorem}\label{Mainequivalence}
		Let $\rho$ be a convex modular on a real vector space $X$, and let $\tau_{\rho}$ and $\tau_{\|\cdot\|_{\rho}}$ be the modular topology and the norm topology, respectively. Then, the following conditions are equivalent:
		\begin{enumerate}
			\item[(i)] $\tau_{\rho}$ is a TVS topology on $X_{\rho}$.
			\item[(ii)] $\rho$ satisfies the $\Delta_2$-property.
			\item[(iii)] $\|\cdot\|_{\rho}$-convergence is equivalent to $\rho$-convergence.
			\item[(iv)] $\tau_{\|\cdot\|_{\rho}} = \tau_{\rho}$.
			\item[(v)] $(X_{\rho}, \tau_{\rho})$ is normable.
		\end{enumerate}
	\end{theorem}
	
	\begin{proof}\leavevmode
		\begin{enumerate}
			\item[(1)] $(i) \Rightarrow (ii)$. Assume $(i)$. Fix a neighborhood of $0$, say $\mathcal{N}$, and let $\mathcal{M}$ be another neighborhood of $0$ such that $\mathcal{M} + \mathcal{M} \subset \mathcal{N}$. By definition, for some $r > 0$ it holds that $B_{\rho, r}(0) \subseteq \mathcal{M}$, and thus
			\[
			2B_{\rho, r}(0) \subset B_{\rho, r}(0) + B_{\rho, r}(0) \subset \mathcal{N}.
			\]
			Pick an arbitrary sequence $(x_j)$ such that $(x_j) \overset{\rho}{\rightarrow} 0$. Then, for some $n_0 \geq 1$, we have $x_j \in B_{\rho, r}(0)$ for $j \geq n_0$, which guarantees $2x_j \in \mathcal{N}$ for $j \geq n_0$. On account of the arbitrariness of $\mathcal{N}$, it follows that $(2x_j)\overset{\tau_{\rho}}{\rightarrow} 0$. By Theorem \ref{convergencecharacterization}, $(2x_j)\overset{\rho}{\rightarrow} 0$. Hence, $\rho$ satisfies the $\Delta_2$-property, as claimed.
			
			\item[(2)] $(ii) \Rightarrow (iii)$. Under the assumption of the $\Delta_2$-property, it follows easily that $(x_j)\overset{\rho}{\rightarrow} x$ implies $(\lambda x_j)\overset{\rho}{\rightarrow} \lambda x$ for any $\lambda > 0$. Proposition \ref{equivalence} yields $\|x_j - x\|_{\rho} \rightarrow 0$ as $j \rightarrow \infty$. Additionally, $\|\cdot\|_{\rho}$-convergence implies $\rho$-convergence. Thus, $(iii)$ holds.
			
			\item[(3)] $(iii) \Rightarrow (iv)$. Recall that, on account of $(iii)$, $\tau_{\rho} \subseteq \tau_{\|\cdot\|_{\rho}}$. Now, let $A$ be $\tau_{\|\cdot\|_{\rho}}$-closed. Then, any sequence $(x_j) \subseteq A$ with $(x_j) \overset{\rho}{\rightarrow} x$ must converge to $x$ in the Luxemburg norm $\|\cdot\|_{\rho}$. The norm-closedness of $A$ implies $x \in A$, so that $A$ is also $\rho$-closed. Thus $\tau_{\|\cdot\|_{\rho}} \subseteq \tau_{\rho}$. Hence, both topologies coincide.
			
			\item[(4)] $(iv) \Rightarrow (v)$. This is immediate.
			
			\item[(5)] $(v) \Rightarrow (i)$. If $(v)$ holds, then $\tau_{\rho}$ is the topology generated by a norm $\|\cdot\|$, which is a TVS topology on $X_{\rho}$.
		\end{enumerate}
		Thus, the proof of Theorem \ref{Mainequivalence} is complete.
	\end{proof}
	
	The following theorem shows that in the finite dimensional case, the modular structure offers no additional insights.
	
	\begin{theorem}\label{equivalence-finte-dim}
		Let $\rho: X \rightarrow [0, \infty]$ be a convex modular on a finite dimensional vector space $X$. Then $\rho$ satisfies the $\Delta_2$-property. Hence, the topology of the Luxemburg norm and the modular topology coincide on any finite dimensional space.
	\end{theorem}
	\begin{proof}
		Let $(x_j) \subset X$ be a sequence for which $\lim\limits_{j \to \infty} \rho(x_j) = 0$. If $(\rho(2x_j))$ did not tend to zero, then for some $\varepsilon_0 > 0$, there would exist a subsequence $(x_{j_k})$ for which $\rho(x_{j_k}) \rightarrow 0$ and $\rho(2x_{j_k}) \geq \varepsilon_0$. Thus, there would exist $k_1$ such that for $k \geq k_1$, it would hold that $\rho(x_{j_k}) < 1$, which implies $\|x_{j_k}\|_{\rho} \leq 1$. Because $X$ is finite dimensional, it has the Heine-Borel property, and this implies that $(x_{j_k})_{k \geq k_1}$ has a norm-convergent subsequence, say $(y_n)$. Since $\rho(y_n) \rightarrow 0$, $\|y_n\|_{\rho}$ must also converge to $0$ as $n \rightarrow \infty$. Thus, $\|2y_n\|_{\rho} \rightarrow 0$ as $n \rightarrow \infty$. This contradicts the fact that $\rho(2y_{n}) \geq \varepsilon_0$ for all $n \in \mathbb{N}$. In all, $\rho(2x_n) \rightarrow 0$, and $\rho$ satisfies the $\Delta_2$-property.
	\end{proof}

	\begin{theorem}\label{equivalence-compact}
		Let $\rho:X\rightarrow [0,\infty]$ be a convex modular on any  vector space $X$ and $K\subset X_\rho$ be compact relative to the $\|\cdot\|_{\rho}$-topology. Then any sequence $(x_j)\subseteq K$ is $\rho$-convergent if and only if it is $\|\cdot\|_{\rho}$-convergent.
	\end{theorem}
	\begin{proof} The proof is straightforward.
	\end{proof}

	Perhaps the most significant shortcoming of the topology $\tau_{\rho}$ is the counterintuitive fact that modular balls are not necessarily open. The following theorem sheds some light on this issue.
	
	\begin{theorem}\label{delta2-rcontinuity}
		Let $\rho$ be a convex, left-continuous modular on a vector space $X$ and consider the following statements:
		\begin{enumerate}
			\item[(i)] each open $\rho$-ball is $\rho$-open;
			\item[(ii)] $\rho$ is right-continuous;
			\item[(iii)] each open $\rho$-ball is norm-open.
		\end{enumerate}
		Then $(i)\Rightarrow (ii)\Rightarrow (iii)$.
		In addition, if $\rho$ satisfies the $\Delta_2$-property, then $(iii)\Rightarrow (i)$.
	\end{theorem}
	
	\begin{proof}\leavevmode
		\begin{enumerate}
			\item[(1)]  $(i) \Rightarrow (ii)$. Assume $(i)$. If $(ii)$ fails, let $x_0\in X_\rho$ be such that for some sequence $(\lambda_j)\searrow 1 $ and $\delta > 0$, we have $\rho(\lambda_j x_0)\geq \rho(x_0)+\delta$ for any $j \in \mathbb{N}$. Consider the open modular ball $B=\{z\in X_\rho:\ \rho(z)<\rho(x_0)+\frac{\delta}{2}\}$. Since $(\lambda_j x_0)\overset{\rho}{\rightarrow}x_0$ and $x_0 \in B$, it is easy to see that $B$ does not contain any $\rho$-ball centered at $x_0$. Hence $B$ is not $\rho$-open, contradicting $(i)$.
			\item[(2)]  $(ii) \Rightarrow (iii)$. Assume $(ii)$ holds. First, establish that for $r>0$ 
			\begin{equation}\label{1tor}
				B_r(0)=\{x \in X_\rho:\ \rho(x)<r\}=\{x \in X_\rho:\ \mu_{B_r}(x)<1\},
			\end{equation}
			where $\mu_{B_r}$ is the Minkowsky functional associated with $B_r(0)$. If $\rho(x)<r$, then by right-continuity, there exists $\lambda>1$ such that $\rho(\lambda x)<r$, i.e., $\mu_{B_r}(x)\leq \lambda^{-1}<1$. Conversely, if 
			$$\mu_{B_r}(x)=\inf\{\gamma>0: \rho\left(\gamma^{-1}x\right)\leq r \}<1,$$
			then for some $\gamma \in (0,1)$, one has $\rho (\gamma^{-1}x)\leq r$, yielding
			$$\rho(x)=\rho(\gamma^{-1}\gamma x)\leq \gamma \rho(\gamma^{-1}x)\leq \gamma r<r.$$
			Thus, (\ref{1tor}) holds. Since the norm $\mu_{B_r}$ is equivalent to the Luxemburg norm, the right-hand side in (\ref{1tor}) is $\|\cdot\|_{\rho}$-open, and $(iii)$ follows.
		\end{enumerate}
		Finally, if the $\Delta_2$ condition holds, then the norm topology and the $\rho$-topology coincide. If that is the case, it is clear that $(iii)\Rightarrow (i)$.
	\end{proof}
	
	A straightforward consequence of Theorem \ref{convergencecharacterization} is that the continuous functions with respect to the modular topology are precisely those that are sequentially continuous. The following result holds:
	
	\begin{theorem}\label{continuitysequential}
		Let $(X,\rho_X)$ and $(Y,\rho_Y)$ be modular spaces with modular topologies $\tau_{\rho_X}$ and $\tau_{\rho_Y}$, respectively. Then
		$$f:(X_{\rho_X},\tau_{\rho_X})\rightarrow (Y_{\rho_Y},\tau_{\rho_Y})$$
		is continuous if and only if it is sequentially continuous.
	\end{theorem}
	\begin{proof}
		Due to Theorem \ref{convergencecharacterization}, modular convergence of sequences is equivalent to modular topological convergence. Therefore, no distinction should be made between the two types of convergence. Assume $f$ is continuous and $(x_j)\overset{\tau_{\rho_X}}{\rightarrow}x$. If $(f(x_j))$ does not converge to $f(x)$ in the $\tau_{\rho_Y}$-topology, then there would be a $\tau_{\rho_Y}$-open neighborhood $U$ of $f(x)$ and a subsequence $(x_{j_k})$ with $(f(x_{j_k}))\subseteq Y\setminus U$. Now, $Y\setminus U$ is $\tau_{\rho_Y}$-closed, and the continuity assumption implies $f^{-1}\left(Y\setminus U\right)$ must be $\tau_{\rho_X}$-closed. Thus, the fact that $(x_{n_k})\subseteq f^{-1}\left(Y\setminus U\right)$ $\rho_X$-converges to $x$ would force $x\in f^{-1}\left(Y\setminus U\right)$, which is impossible. It follows that the sequence $(f(x_j))$ must converge to $f(x)$ in the $\tau_{\rho_Y}$-topology, as claimed.\\
		Conversely, assume $f$ is sequentially continuous. Let $C\subset Y$ be $\tau_{\rho_Y}$-closed. We will prove that $f^{-1}(C)$ is $\tau_{\rho_X}$-closed. Let $(x_j)\subseteq f^{-1}(C)$ with $(x_j)\overset{\tau_{\rho_X}}{\rightarrow}x$. Under the assumption of sequential continuity, it must hold $f(x_j)\overset{\rho_{\tau_{\rho_Y}}}{\rightarrow}f(x).$ Since $C$ is closed, it follows $f(x)\in C$, i.e., $x\in f^{-1}(C)$. This implies $f^{-1}(C)$ is $\tau_{\rho_X}$-closed, as claimed.
	\end{proof}
	\section{The modular topology of $\ell^{(p_n)}$ and $L^{p(\cdot)}$}\label{lpn}
	In this section the modular topologies of the variable exponent sequence spaces and Lebesgue spaces (Definitions \ref{deflp} and \ref{Lp}) are studied in detail.  
	
	\begin{theorem}\label{r-continuous-p^+}
		Let \( \mathbf{p}: = (p_j) \subset (1,\infty) \) be a sequence. Recall that the modular \(\rho_{\mathbf{p}}: \ell^{(p_j)} \to [0,\infty]\) is defined as follows:
		$$\rho_{\mathbf{p}}((a_j)):=\textstyle \sum\limits_{j=1}^{\infty}|a_j|^{p_j}.$$
		Then the following are equivalent:
		\begin{enumerate}
			\item[(i)] $\rho_{\mathbf{p}}$ is right-continuous on $\ell^{(p_j)}$,
			\item[(ii)] $p_+:=\sup\limits_{j\in {\mathbb N}}\ p_j<\infty$,
			\item[(iii)] $\rho_{\mathbf{p}}$ satisfies the $\Delta_2$-property.
		\end{enumerate}
	\end{theorem}
	\begin{proof}
		We will first show that if $p_+ = \infty$, then $\rho_{\mathbf{p}}$ is not right-continuous. There exists a strictly increasing sequence \((n_k)\subseteq \mathbb{N}\), with \(n_k\geq k\), such that \(p_{n_k}>k^2\) and that \((p_{n_k})\) is strictly increasing. Define the sequence \((a_k)\) by setting \(a_{n_k} = p_{n_k}^{-\frac{1}{p_{n_k}}}\) for \(k=1,2,\ldots\) and \(a_j = 0\) for \(j\neq n_k\). Then \((a_j) \in \ell^{(p_n)}\), as 
		$$\rho_{\mathbf{p}}((a_j)) = \textstyle \sum\limits_{k=1}^{\infty} \frac{1}{p_{n_k}} \leq \textstyle \sum\limits_{k=1}^{\infty} \frac{1}{k^2} < \infty.$$
		On the other hand, for any \(\lambda > 1\), we have
		$$\rho_{\mathbf{p}}((\lambda a_j)) = \textstyle \sum\limits_{k=1}^{\infty} \frac{\lambda^{p_{n_k}}}{p_{n_k}} = \infty.$$
		It follows that \(\rho_{\mathbf{p}}\) is not right-continuous, thus \((i) \Rightarrow (ii)\).\\
		Conversely, assume \((ii)\) holds, i.e., \(p_+ < \infty\). For any \((a_j) \in \ell^{(p_n)}\), the following inequality holds:
		$$\rho_{\mathbf{p}}((2a_j)) \leq 2^{p_+} \rho_{\mathbf{p}}((a_j)).$$
		If \(\rho_{\mathbf{p}}((a_j)) = \infty\), then for any \(\lambda > 1\), it is clear that
		$$\rho_{\mathbf{p}}(\lambda (a_j)) \geq \rho_{\mathbf{p}}((a_j)) = \infty.$$
		Furthermore, it is clear that \(\lim_{\lambda \to 1^+} \rho_{\mathbf{p}}(\lambda (a_j)) = \rho_{\mathbf{p}}((a_j))\). If \(\rho_{\mathbf{p}}((a_j)) < \infty\) and \(\lambda_k \searrow 1\), then, for some \(\delta > 0\) and for each \(k\), \((\lambda_k)^{p_j} < (1 + \delta)^{p_+}\). Thus, for any \(j \geq 1\),
		$$|\lambda_k a_j|^{p_j} \leq (1 + \delta)^{p_+} |a_j|^{p_j},$$
		and the series is summable. On account of Lebesgue's dominated convergence theorem it follows that
		$$\lim_{\lambda \to 1^+} \textstyle \sum\limits_{j=1}^{\infty} |\lambda_k a_j|^{p_j} = \textstyle \sum\limits_{j=1}^{\infty} |a_j|^{p_j}.$$
		Thus, \((i) \iff (ii)\). The implication \((ii) \Rightarrow (iii)\) follows directly. Finally, if \((ii)\) is not satisfied, then \((iii)\) cannot hold. If the sequence \((p_n)\) is unbounded, a strictly increasing subsequence \((n_k)\) of natural numbers can be chosen such that \(p_{n_k} > k\) and \((p_{n_k})\) is strictly increasing. Let \((x_k)\) be the sequence equal to \(\frac{1}{2}\) on each \(n_k\) and zero otherwise. Then the sequence \(\Big((x_k)\mathbbm{1}_{\{i:\ i > n_m\}}\Big)_{m \geq 1}\) \(\rho_{\mathbf{p}}\)-converges to \(0\), but \(\Big(2 (x_k)\mathbbm{1}_{\{i:\ i > n_m\}}\Big)_{m \geq 1}\) does not. Hence, \((iii)\) fails, proving that \((iii) \Rightarrow (ii)\).
	\end{proof}
	A similar situation arises in the continuous case:
	
	\begin{theorem}\label{r-continuous-p^+L}
		Let $L^{p(\cdot)}(\Omega)$ and $\rho_{{p}}$ be as in Definition \ref{Lp}.
		Then the following are equivalent:
		\begin{enumerate}
			\item[(i)] $\rho_{{p}}$ is right-continuous,
			\item[(ii)] $p_+=\|p\|_{\infty}<\infty$,
			\item[(iii)] $\rho_{{p}}$ satisfies the $\Delta_2$-property.
		\end{enumerate}
	\end{theorem}
	\begin{proof}
		The theorem follows along the same lines as the preceding one and will only be sketched. For each $k\in {\mathbb N}$, $\Omega_k=\{x\in \Omega :k\leq p(x)<k+1\}$. If $(ii)$ fails, then $|\Omega_{k_j}|\neq 0$ for an infinite sequence $\{k_1,...k_j,...\}$ and the function $u=\sum\limits_{j=1}^{\infty}{\mathbbm 1}_{\Omega_{k_j}}2^{-1}|\Omega_{k_j}|^{-\frac{1}{p(x)}}$ belongs to $L^{p(\cdot)}(\Omega)$, $\rho_{{p}}(u)<1$ and $\rho_{{p}}(2u)=\infty.$ Thus, $(iii)\Rightarrow (ii)$. Next, for 
		$$v(x)=\sum\limits_{j=1}^{\infty}{\mathbbm 1}_{\Omega_{k_j}}n^{-2}|\Omega_{k_j}|^{-\frac{1}{p(x)}}$$
		it is clear that $\rho_{{p}}(\theta v)=\infty$ for any $\theta>1$,, whereas $\rho_{{p}}(v)<\infty.$ Thus, $(i)\Rightarrow (ii)$. The implication $(iii)\Rightarrow (i)$ follows from Dominated Convergence. The proof of $(ii)\Rightarrow (iii)$ is obvious.
	\end{proof}
	
	The following corollary is a direct consequence of Theorems  \ref{delta2-rcontinuity} and \ref{equivalence}:
	
	\begin{corollary}
		In the notation of Theorem \ref{r-continuous-p^+}, all (open) modular balls in \(\ell^{(p_j)}\) ($L^{p(\cdot)}(\Omega)$) are \(\tau_{\rho_{\mathbf{p}}}\)-open if and only if the sequence \((p_j)\) (the exponent function $p(x)$) is bounded, i.e., \(p_+ = \sup\limits_j p_j < \infty\) ($p_+=\|p\|_{\infty}<\infty$).
	\end{corollary}
	
	The next two examples strengthen the previous result: they show that that in the absence of the \(\Delta_2\)-property, no open modular ball is \(\rho_{\mathbf{p}}\)-open in $\ell^{(p_n)}$.
	
	\begin{example}\label{nice-example}
		{\normalfont 
			Consider the variable exponent space \(\ell^{(p_n)}\), where \(p_n = n\) for all \(n \geq 1\). The sequence \((x_n) \subset \ell^{(p_n)}\) is defined such that \(x_n\) takes the value \(\frac{1}{2}\) up to the \(n^{th}\)-position and \(0\) for all other positions. For \(n \geq 1\), we have
			$$\rho(x_n) = \textstyle \sum\limits_{j=1}^{n}\left(\frac{1}{2}\right)^j.$$
			Now, consider the modular ball of radius 3 centered at \({\mathbbm 1}_{{\mathbb N}}\), which is given by 
			$$B_{\rho,3}({\mathbbm 1}_{{\mathbb N}}) = \Big\{ (z_j) \in \ell^{(p_n)}:\ \rho((z_j-1)) < 3\Big\}.$$ 
			This ball contains the sequence \(x = \frac{1}{2} {\mathbbm 1}_{{\mathbb N}}\). Since \((x_n)\) $\rho$-converges to \(x\), any modular ball centered at \(x\) must contain some element \(x_n\), for sufficiently large \(n\), which contradicts the fact that none of the \(x_n\)'s are in \(B_{\rho,3}({\mathbbm 1}_{{\mathbb N}})\). Hence, the ball \(B_{\rho,3}({\mathbbm 1}_{{\mathbb N}})\) is not $\tau_{\rho}$-open.\\
			More generally, for any \(\varepsilon > 0\), no \(x = (x_j) \in B_{\rho,\varepsilon}({\mathbbm 1}_{{\mathbb N}})\) with \(\textstyle \sum\limits_{1}^{\infty}|x_j|^j < \infty\) is an interior point of \(B_{\rho,\varepsilon}({\mathbbm 1}_{{\mathbb N}})\). To see this, take any \(\delta > 0\). The sequence \(y_n = (x_1, x_2, \dots, x_n, 0, 0, \dots)\) will belong to \(B_{\rho,\delta}(x)\) for sufficiently large \(n\), but clearly, \(y_n \notin B_{\rho,\varepsilon}({\mathbbm 1}_{{\mathbb N}})\). Thus, no modular ball is \(\rho\)-open. In fact, for any \(x \in \ell^{(p_n)}\) and \(\varepsilon > 0\),
			$$B_{\rho,\varepsilon}(x) = \big(x - {\mathbbm 1}_{{\mathbb N}}\big) + B_{\rho,\varepsilon}({\mathbbm 1}_{{\mathbb N}}),$$
			and the result follows from (1) in Proposition \ref{basic-properties}.}
	\end{example}
	
	The conclusion from Example \ref{nice-example} extends to any space \(\ell^{(p_n)}\) when \(p_n \to \infty\).
	
	\begin{example}\label{modularballsnotopen}
		{\normalfont
			For arbitrary \((p_n)\), where \(p_n \to \infty\), construct the sequence \((n_k)\) as follows: choose \(p_{n_1} \geq 1\) and \(p_{n_k} \geq \max\{k, p_{n_{k-1}}\}\) for \(k > 1\). The sequence \(\mathbf{s} = {\mathbbm 1}_{\{n_1, n_2, ..., n_k, ...\}}\) belongs to \(\ell^{(p_n)}\) because
			$$\textstyle \rho_{\mathbf{p}}\left(2^{-1} \mathbf{s}\right) = \textstyle \sum\limits_1^{\infty}2^{-p_{n_k}} \leq \textstyle \sum\limits_1^{\infty}\ 2^{-k} < \infty.$$
			For any \(\delta > 0\), \((1 - \varepsilon) \mathbf{s} \in B_{\rho_{\mathbf{p}}, \delta}(\mathbf{s})\) for sufficiently small \(\varepsilon\). Notice also that 
			$$\rho_{\mathbf{p}}\left((1 - \varepsilon) \mathbf{s}\right) < \infty.$$
			Thus, \((1 - \varepsilon) \mathbf{s}\) can be approximated in the modular sense by
			$$\left((1 - \varepsilon)(\mathbf{s}_1, \dots, \mathbf{s}_N, 0, \dots)\right),$$ 
			which is \textbf{not} in the ball \(B_{\rho_{\mathbf{p}},\delta}(\mathbf{s})\). Therefore, the modular ball \(B_{\rho_{\mathbf{p}},\delta}(\mathbf{s})\) is not open in the modular topology \(\tau_{\rho_{\mathbf{p}}}\). As before, no modular ball is open in the modular topology of the \(\ell^{(p_n)}\) space.}
	\end{example}
	\begin{definition}
		Given a modular space $(X,\rho)$ and a subset $Y\subseteq A$, the modular diameter of $Y$, is defined as $$\text{diam}_{\rho}(Y):=\sup\{\rho(a-b),a\in Y, b\in Y\}$$.
		
	\end{definition}
	\begin{remark}{\normalfont 
			A direct consequence of Example \ref{modularballsnotopen} is that non-empty open sets in the modular topology must have infinite modular diameter in \(\ell^{(p_n)}\). This is because, as shown in the example, for any \(\delta > 0\), the modular ball \(B_{\rho_{\mathbf{p}},\varepsilon}(\mathbf{s})\) contains a point \(\mathbf{t}\) such that \(\rho_{\mathbf{p}}(\mathbf{s} - \mathbf{t}) = \infty\). Therefore, if \(\mathcal O\) is modularly open and \(w \in \mathcal{O}\), the set \(\mathcal{O} + \{\mathbf{s} - w\}\) is modularly open and contains \(\mathbf{s}\), and hence, a modular ball \(B_{\rho_{\mathbf{p}},\delta}(\mathbf{s})\). It follows that
			$$\text{diam}_\rho\left(\mathcal{O} + \{\mathbf{s} - w\}\right) = \infty,$$
			and since the \(\rho\)-diameter is translation-invariant, \(\mathcal{O}\) must have infinite diameter. Recall that for any non-empty subset \(A\) in \(\ell^{(p_n)}\), the \(\rho\)-diameter of \(A\) is defined as
			$$\text{diam}_\rho(A) := \sup\limits_{x,y \in A} \rho(x - y).$$}
	\end{remark}
	The following is an extension of Example \ref{nice-example} to the case of the \(L^{p(\cdot)}\) spaces.
	
	\begin{example}\label{nice-example1}{\normalfont  
			Let \(\Omega = (0,1)\), \(p(x) = x^{-1}\), and \(\rho_{p}(v) = \bigint_0^1 |v(x)|^{p(x)} dx\). Define 
			$$v(x): = \textstyle \sum\limits_{n=1}^{\infty} n^{\frac{1}{n}} \cdot {\mathbbm 1}_{((n+1)^{-1}, n^{-1})}.$$
			It is straightforward to verify that \(\rho_{{p}}(v) = \infty\) and that \(\rho_{{p}}((1 - \varepsilon) v) \to 0\) as \(\varepsilon \to 1\). In particular, for \(n \geq 1\), we have
			$$n^{(xn)^{-1}} \leq n \quad \text{and} \quad (1 - \varepsilon)^{x^{-1}} \leq (1 - \varepsilon)^n,$$
			for any \(x \in ((n+1)^{-1}, n^{-1})\), which implies
			\begin{align*}
				\int_{(n+1)^{-1}}^{n^{-1}} n^{(xn)^{-1}} \  (1 - \varepsilon)^{x^{-1}} dx &\leq n (1 - \varepsilon)^n \left(\frac{1}{n} - \frac{1}{n+1}\right) = \frac{(1 - \varepsilon)^n}{n + 1} \\
				&< (1 - \varepsilon)^n.
			\end{align*}
			Hence,
			$$\rho_{{p}}((1 - \varepsilon) v) < \textstyle \sum\limits_{n=1}^{\infty} (1 - \varepsilon)^n = \displaystyle  \frac{1 - \varepsilon}{\varepsilon},$$
			which proves the claim. On the other hand, for \(k \geq 1\), define
			$$v_k := {\mathbbm 1}_{((k+1)^{-1}, 1)} v = \textstyle \sum\limits_{n=1}^{k} n^{n^{-1}} \cdot {\mathbbm 1}_{((n+1)^{-1},n^{-1})}.$$
			For any fixed \(k \geq 1\) and \(\varepsilon \in (0,1)\), we have \(\rho_{{p}}(v - v_k) = \infty\). Using similar arguments as before, it is readily obtained that
			$$\begin{array}{lll}   
				\rho_{{p}}(\varepsilon (v-v_k)) & \displaystyle =\rho\Big( \textstyle \sum\limits_{n=k+1}^ {\infty}\varepsilon \ n^{n^{-1}}\ {\mathbbm 1}_{((n+1)^{-1},n^{-1})}\Big)\\
				&\leq  \textstyle \sum\limits_{n=k+1}^ {\infty}\displaystyle  \frac{\varepsilon^n n}{n(n+1)}\\
				&\displaystyle < \textstyle \sum\limits_{n=k+1}^ {\infty}\varepsilon^n =\displaystyle  \frac{\varepsilon^{k+1}}{1-\varepsilon},
			\end{array}$$
			which implies \(\rho_{{p}}\left(\varepsilon (v - v_k)\right) \to 0\) as \(k \to \infty\). Moreover, \(\rho_{{p}}(v - \varepsilon v_k) = \infty\) for any \(k \geq 1\) and \(\varepsilon \in (0,1)\).  
			Thus, for any \(\delta > 0\), the modular ball \(B_{p, \delta}(v)\) contains \(\varepsilon v\) for some \(\varepsilon \in (0,1)\), yet any modular ball \(B_{p, \alpha}(\varepsilon v)\) contains an element not in \(B_{p, \delta}(v)\). Therefore, \(B_{p, \delta}(v)\) is not \(\rho_{{p}}\)-open.}
	\end{example}
	A stronger result holds: Open modular balls in $\ell^{(p_n)}$ and $L^{p(\cdot)}$ have empty interior when $p_+=\infty$. This follows from Theorem \ref{separability}.
	\begin{theorem}\label{separability}
		$(\ell^{(p_n)},\tau_{\rho_{\mathbf{p}}})$ and $(L^{p(\cdot)}(\Omega),\tau_{\rho_{{p}}})$ are $\tau_{\rho}$- separable. Moreover, $C_0^{\infty} (\Omega)$ is $\tau_{\rho}$-dense in $L^{p(\cdot)}(\Omega)$ if $1\leq p<\infty$, even if $p_+=\infty$.
	\end{theorem}
	\begin{proof}
		Let $S$ be the subspace of $\ell^{(p_n)}$ consisting of sequences that have finitely many non-zero terms and $\overline{S}^\rho$ be its $\rho_{p}$ closure, which is a $\rho_{p}$-closed subspace of $\ell^{(p_n)}$.  If $\rho_{{p}}((a_j))<\infty$, it is plain that any modular ball $B_{\delta}((a_j))$ contains an element of $S$. By definition, any $\tau_{\rho_{{p}}}$ neighborhood of $(a_j)$ contains a modular ball, from which it follows that any $\tau_{\rho}$-open set containing $(a_j)$ has nonempty intersection with $S$. In other words, $(a_j)\in \overline{S}^\rho$. Next, observe that for any $(b_j)\in \ell^{(p_n)}$ there exists $\lambda>0$ such that $\rho_{p}(\lambda (b_j))<\infty,$ i.e., $(b_j)\in \overline{S}^\rho$.\\
		For the remaining part of the proof, let $\Omega \subseteq {\mathbb R}^n$ be a domain and $C_0^{\infty}(\Omega)$ the vector space of compactly supported, infinitely differentiable functions on $\Omega$.\\ It is apparent that it is enough to show that given any measurable $f$ on $\Omega$ with $\int\limits_{\Omega}|f(x)|^{p(x)}dx<\infty$ for any $\delta>0$ it holds
		\begin{equation}\label{sep}
			B_{\delta}(f)\cap C_0^{\infty}(\Omega)\neq \emptyset.
		\end{equation}
		Let $\Omega^{\prime}\subset \Omega$ be chosen so that $\int\limits_{\Omega^{\prime}}|f(x)|^{p(x)}dx<\delta/2.$ Clearly, $f\in L^{p{\mathbbm 1}_{\Omega\setminus \Omega^{\prime}}}(\Omega\setminus \Omega^{\prime})$. Since $C^{\infty}_0(\Omega\setminus \Omega^{\prime})$ is dense in $L^{p{\mathbbm 1}_{\Omega\setminus \Omega^{\prime}}}(\Omega\setminus \Omega^{\prime})$, there exists $\phi \in C^{\infty}_0(\Omega\setminus \Omega^{\prime})$ such that
		\begin{equation}
			\int\limits_{\Omega\setminus \Omega^{\prime}}|f(x)-\phi(x)|^{p(x){\mathbbm 1}_{\Omega\setminus \Omega^{\prime}}}dx<\delta/2.
		\end{equation}
		Then,
		\begin{align*}
			\int\limits_{\Omega}|f(x)-\phi(x)|^{p(x)}dx&=\int\limits_{\Omega^{\prime}}|f(x)-\phi(x)|^{p(x)}dx+\int\limits_{\Omega\setminus \Omega^{\prime}}|f(x)-\phi(x)|^{p(x)}dx \\ & =\int\limits_{\Omega^{\prime}}|f(x)|^{p(x)}dx +\delta/2\\
			& <\delta.
		\end{align*}
		This concludes the proof.
	\end{proof}
	
	\begin{corollary}
		If the sequence ${\mathbf{p}}=(p_n)$ is unbounded, then open modular balls in $\ell^{(p_n)}$ have empty interior.
	\end{corollary}
	\begin{proof}
		It was shown in Example \ref{nice-example} that for any $\delta>0$, the open modular ball $B_{\delta}({\mathbbm 1}_{\mathbb N})$ in $\ell^{(p_n)}$ with $(p_n)=(n)$, has no interior point $\mathbf {x}$ with $\rho_{\mathbf{p}}(\mathbf{x})<\infty.$ Select $\mathbf {y}\in B_{\delta}({\mathbbm 1}_{\mathbb N})$. On account of Theorem \ref{separability}, any $\rho_{\mathbf{p}}$-open set $A$ containing $\mathbf {y}$ must contain a sequence with only finitely many nonzero terms. Pick one such sequence, say $\mathbf{x}$; then  $\mathbf{x}$ is not an interior point of $B_{\delta}({\mathbbm 1}_{\mathbb N})$, i.e., some point of $A$ must be in the complement of $B_{\delta}({\mathbbm 1}_{\mathbb N})$.
	\end{proof}

	\medskip
	\subsection{Separation properties}\leavevmode
	
	We point out that, by definition, every singleton \(\{a\} \subset X_{\rho}\) is \(\rho\)-closed for any modular topology. 
	
	\begin{definition}
		The modular \(\rho\) on the vector space \(X\) is said to satisfy the Fatou property if whenever \((y_j) \overset{\rho}{\rightarrow} y \in X_{\rho}\), it holds that \(\rho(y) \leq \liminf\limits_{j \to \infty} \rho(y_j).\)
	\end{definition}
	
	It is evident that the Fatou property holds if and only if the modular balls \(\{v \in X_\rho : \rho(x - v) \leq \varepsilon\}\), for any \(\varepsilon > 0\) and \(x \in X_\rho\), are closed in \(\tau_{\rho}\).
	
	\begin{lemma}\label{t0}
		If \(\rho\) satisfies the Fatou property, then the topology \(\tau_{\rho}\) is \(T_1\).
	\end{lemma}
	\begin{proof}
		Let \(a\) and \(b\) be distinct points in \(X_{\rho}\), i.e., \(a \neq b\). Set
		$$\mathcal{O} = \left\{x : \rho(x - a) > \rho\left(\frac{a - b}{2}\right)\right\}.$$
		On account of the Fatou property, \(\mathcal{O}\) is \(\tau_{\rho}\)-open. Clearly, \(a \notin \mathcal{O}\) and \(b \in \mathcal{O}\).
	\end{proof}
	
	\begin{question}
		Is \(\tau_{\rho}\) regular (i.e., \(T_3\))?
	\end{question}
	
	\begin{example}\label{lphausdorff}
		{\normalfont 
			Consider the variable exponent space \(\ell^{(p_n)}\) and the modular \(\rho_{\mathbf{p}}\). By virtue of Lemma \ref{t0}, \(\tau_{\rho_{\mathbf{p}}}\) is \(T_1\). Fix \(M \in \mathbb{N}\). For a given \(x = (x_n) \in \ell^{(p_n)}\) and \(\varepsilon > 0\), the set 
			$$U^M_{x,\varepsilon} = \Big\{y = (y_n) \in \ell^{(p_n)} : |x_M - y_M| < \varepsilon\Big\}$$
			is \(\tau_{\rho_{\mathbf{p}}}\)-open and contains \(x\). To demonstrate this, let us prove that \(\ell^{(p_n)} \setminus U^M_{x,\varepsilon}\) is \(\tau_{\rho_{\mathbf{p}}}\)-closed. Consider a sequence \((\psi^n) \subset \ell^{(p_n)} \setminus U^M_{x,\varepsilon}\) and assume that \((\psi^n) \overset{\rho_{\mathbf{p}}}{\rightarrow} \psi.\) Then, for any \(\delta > 0\), there exists \(N \geq 1\) such that
			$$\delta^{p_M} > \rho(\psi^n - \psi) = \sum_{j=1}^{\infty} |\psi_j^n - \psi_j|^{p_j} \geq |\psi^n_M - \psi_M|^{p_M},$$
			for \(n > N\). Thus, for \(n > N\),
			$$|\psi_M - x_M| \geq |\psi^n_M - x_M| - |\psi^n_M - \psi_M| \geq \varepsilon - \delta.$$
			Since \(\delta\) was chosen arbitrarily, we have \( |\psi_M - x_M| \geq \varepsilon \), i.e., \(\psi \in \ell^{(p_n)} \setminus U^M_{x,\varepsilon}\), which proves the claim. It is now a routine matter to verify that \(\tau_{\rho}\) is Hausdorff (\(T_2\)). Indeed, if \(a\) and \(b\) are in \(\ell^{(p_n)}\), with \(a \neq b\), then \(a_M \neq b_M\) for some \(M \in \mathbb{N}\). The sets \(U^M_{a,\varepsilon}\) and \(U^M_{b,\varepsilon}\), with \(\varepsilon = \frac{|a_M - b_M|}{2}\), are disjoint \(\tau_{\rho_{\mathbf{p}}}\)-open neighborhoods of \(a\) and \(b\), respectively.}
	\end{example}
	
	The only noteworthy case in Example \ref{lphausdorff} arises when the sequence \((p_n)\) is unbounded. For a bounded exponent sequence \((p_n)\), it is straightforward to show that the modular \(\rho\) satisfies the \(\Delta_2\)-property, and according to Theorem \ref{Mainequivalence}, \(\tau_{\rho}\) corresponds to the topology induced by the Luxemburg norm. Therefore, it is \(T_i\) for \(i = 0, 1, 2, 3, 4\).
	
	\begin{lemma}\label{lemmaslow}
		Let \((p_n) \subset (1,\infty)\) be unbounded. Then:
		\begin{itemize}
			\item[$(i)$] The inclusion
			$$(\ell^{(p_n)},\|\cdot\|_{(p_n)}) \hookrightarrow (\ell^{\infty},\|\cdot\|_{\infty})$$ 
			is continuous \cite{orlicz1931}.
			\item[$(ii)$] The inclusion
			$$(\ell^{(p_n)},\|\cdot\|_{(p_n)}) \hookrightarrow (\ell^{\infty},\|\cdot\|_{\infty})$$ 
			is an isomorphism if and only if there exists a constant \(\lambda \in (0,1)\) such that 
			\begin{equation}\label{condinfty}
				\textstyle \sum\limits_{n=1}^{\infty} \lambda^{p_n} < \infty.
			\end{equation}
			\item[$(iii)$] If \(\ell^{(p_n)} \subsetneq \ell^{\infty}\), there exists a proper subsequence \((p_{n_k})\), say $(q_k)$, such that the inclusion
			\(i_{q,\infty}:\left(\ell^{(q_k)}, \|\cdot\|_{q(\cdot)}\right)\rightarrow \left(\ell^{\infty},\|\cdot\|_{\infty}\right)\) is an isomorphism.
			\item[$(iv)$] There exists an infinite, proper subsequence \((p^*_{n_k})=(q_k)\) such that the inclusion  \(i_{p,\infty}:\left(\ell^{(q_k)},\|\cdot\|_{q(\cdot)}\right) \rightarrow \left(\ell^{\infty},\|\cdot\|_{\infty}\right)\)
			is continuous but not onto.
		\end{itemize}
	\end{lemma}
	\begin{proof}\leavevmode
		\begin{itemize}
			\item[\((i)\)] The set-theoretic inclusion \(\ell^{(p_n)} \subseteq \ell^{\infty}\) is clear: if \((a_k) \in \ell^{(p_n)}\), there must exist \(\lambda > 0\) such that 
			$$\textstyle \sum\limits_{n=1}^\infty |\lambda a_n|^{p_n} < \infty.$$
			Hence, \((\lambda a_n)\) is bounded, and since \(\lambda > 0\), it follows that \((a_k) \in \ell^{\infty}\). Moreover, it is also easily seen from the above reasoning that if 
			$$\textstyle \sum\limits_{n=1}^{\infty} |\lambda a_n|^{p_n} \leq 1,$$ 
			then \(|a_n| \leq \lambda^{-1}\) for all \(n \in \mathbb{N}\), that is,
			$$\sup\limits_{n \in \mathbb{N}} |a_n| = \|(a_n)\|_{\infty} \leq \lambda^{-1}.$$
			By definition, it follows that $\|(a_n)\|_{\infty} \leq \|(a_n)\|_{(p_n)}$.
			\item[\((ii)\)] If the inclusion in \((i)\) is an isomorphism, it must hold that 
			\({\mathbbm 1}_{{\mathbb N}} \in \ell^{(p_n)}\); hence, by virtue of Definition \ref{deflp}, condition (\ref{condinfty}) must hold. Conversely, under the assumption \((\ref{condinfty})\), observe that if \((b_n) \in \ell^{\infty}\), then
			\begin{equation}\label{lambdabn}
				\textstyle \sum\limits_{n=1}^{\infty} \left|\frac{\lambda}{\|(b_n)\|_{\infty}} b_n\right|^{p_n} \leq \textstyle \sum\limits_{1}^{\infty} \lambda^{p_n}.
			\end{equation}
			In other words, \((b_n) \in \ell^{(p_n)}\). If \(\textstyle \sum\limits_{1}^{\infty} \lambda^{p_n} \leq 1\), inequality \((\ref{lambdabn})\) yields 
			\begin{equation}\label{pninfty}
				\|(b_n)\|_{(p_n)} \leq \lambda^{-1} \|(b_n)\|_{\infty}.
			\end{equation}
			On the other hand, \((\ref{pninfty})\) also holds if \(\textstyle \sum\limits_{1}^{\infty} \lambda^{p_n} > 1\), for in this case, one has
			$$\textstyle \sum\limits_{n=1}^{\infty} \left|\frac{\lambda b_n}{\|(b_n)\|_{\infty} \textstyle \sum\limits_{1}^{\infty} \lambda^{p_j}}\right|^{p_n} \leq \textstyle \sum\limits_{n=1}^{\infty} \left|\frac{\lambda b_n}{\|(b_n)\|_{\infty} \left(\textstyle \sum\limits_{1}^{\infty} \lambda^{p_j}\right)^{\frac{1}{p_n}}}\right|^{p_n} \leq 1.$$
			This completes the proof.
			\item[$(iii)$] Since \((p_n)\) is unbounded, there exists a proper subsequence \((p_{n_k})\) with \(p_{n_k} > k\). For any \(\lambda : 0 < \lambda < 1\), 
			$$\textstyle \sum\limits_{k=1}^{\infty} \lambda^{p_{n_k}} < \infty,$$
			and on account of \((ii)\), it follows that \(\ell^{p_{n_k}} \equiv \ell^{\infty}.\)
			\item[$(iv)$]  It is clear that the subsequence \((p^*_{n_k})\) resulting from removing the elements of the subsequence \((p_{n_k})\) from \((p_n)\) will give \(\ell^{(p_{n_k})} \subsetneq \ell^{\infty}\) since \(\ell^{(p_n)} \subsetneq \ell^{\infty}\). Observe that since \(\ell^{(p_{n_k})} = \ell^{\infty}\), we have \({\mathbbm 1}_{\{p_{n_1}, p_{n_2}, \dots \}} \in \ell^{(p_n)}\).
		\end{itemize}
	\end{proof}
	
	\begin{corollary}
		For \((p_n) \subseteq (1, \infty)\), \((\ell^{(p_n)}, \|\cdot\|_{p_n})\) is isomorphic to \((\ell^{\infty}, \|\cdot\|_{\infty})\) if and only if \({\mathbbm 1}_{{\mathbb N}} \in \ell^{(p_n)}\).
	\end{corollary}
	
	Recall that if \((X,\tau)\) is a first-countable topological space and \(A \subset X\), then for any \(x\) in the \(\tau\)-closure of \(A\), there exists a sequence \((a_j) \subset A\) that \(\tau\)-converges to \(x\).
	
	\begin{lemma}
		In general, $(\ell^{(p_n)},\tau_{\rho_{\mathbf{p}}})$ and $(L^{p(\cdot)}(\Omega),\tau_{\rho_{{p}}})$ are not first-countable.
	\end{lemma}
	\begin{proof}
		Let us take $p_n = n$, for any $n \in \mathbb{N}$.  Consider the subspace $A \subset \ell^{(p_n)}$ consisting of sequences $(a_n)$ with only a finite number of non-zero terms. By Proposition \ref{subspace}, the closure $\overline{A}^{\rho}$ is also a subspace. The sequence $(x_n) = \frac{1}{2}{\mathbbm 1}_{{\mathbb N}}$ belongs to $\overline{A}^{\rho}$ since $\Big(\frac{1}{2}{\mathbbm 1}_{\{j \in {\mathbb N}, j \leq n\}}\Big)$ $\rho$-converges to $(x_n)$. It follows that $2(x_n) = {\mathbbm 1}_{\mathbb N} \in \overline{A}^{\rho}$, but clearly ${\mathbbm 1}_{\mathbb N}$ is not the limit of any sequence in $A$. Therefore, $\ell^{(p_n)}$ is not first-countable. The proof for $L^{p(\cdot)}(\Omega)$ is similar.
	\end{proof}
	
	Since second countability implies first countability, it is immediate that:
	
	\begin{corollary}
		$(\ell^{(p_n)},\tau_{\rho_{\mathbf{p}}})$ and $(L^{p(\cdot)}(\Omega),\tau_{\rho_{{p}}})$ are not second countable.
	\end{corollary}

	\begin{question}
		Is $\ell^{(p_n)}$ regular? (Notice this question should be easier than question 1.)
	\end{question}
	
	\medskip

	\section{Duality in Modular Vector Spaces} \label{duality}
	In this section we focus on  generalizing the concept of a "dual space" with respect to modular topology on \((X_{\rho}, \tau_{\rho})\).

	Let \(\rho\) be a convex modular on a vector space \(X\), and \(\|\cdot\|_{\rho}\) be the Luxemburg norm on \(X_{\rho}\). As usual, the dual of the normed space $(X_\rho,\|\cdot\|_{\rho})$ will be denoted by $X^{\ast}$. Denote with $X^{\rho}$ the collection of all linear functions on $X$ that are continuous with respect to the modular topology $\tau_{\rho}$. The following theorem holds:
	
	\begin{theorem}\label{modular-duality}
		Using the notation from the previous paragraph, \(X^{\rho}\) is a vector space that is generally strictly contained within \(X^{\ast}\). Let \(\Lambda\) be an algebraic linear functional on the modular space \((X_{\rho}, \rho)\). Additionally, consider the following statements:
		\begin{enumerate}
			\item[(i)] $\Lambda \in X^{\rho}$.
			\item[(ii)] $\Lambda$ is continuous at $0$ with respect to the modular topology.
			\item[(iii)] For any sequence $(x_j)\in X$, $(x_j)\overset{\rho}{\rightarrow}0$ implies $\Lambda(x_j)\rightarrow 0$.
			\item[(iv)] $\Lambda$ is bounded on $\{x:\rho(x)<1\}$.
			\item[(v)] $\Lambda$ is bounded on $\{x:\rho(x)\leq 1\}$.
			\item [(vi)] $\Lambda$ is bounded on every modular ball $\{x:\rho(x-x_0)\leq r\}$.
			\item [(vii)] $\Lambda$ is bounded on every modular ball $\{x:\rho(x-x_0)< r\}$.
			\item[(viii)] $\Lambda$ is bounded on every set of finite modular diameter.
			\item [(ix)] $\Lambda \in X^{\ast}$.
		\end{enumerate}
		Then $(i)\iff (ii)\iff (iii)\Rightarrow (iv)\iff (v)\iff (vi)\iff (vii)\Rightarrow (viii)$ and $(i)\Rightarrow (ix)$.
	\end{theorem}
	
	\begin{proof}
		It is obvious that $(i)\Rightarrow (ii)$, and the implication $(ii)\Rightarrow (iii)$ follows from Theorem \ref{continuitysequential}, by observing that the usual topology in ${\mathbb R}$ is the modular topology corresponding to the modular $\eta(s)=|s|$. To show that $(iii)\Rightarrow (i)$, observe that for $X \supset (x_j)\overset{\rho}{\rightarrow x}$, we have $x_j - x \overset{\rho}{\rightarrow}0$. Then $\Lambda(x_j - x) = \Lambda(x_j) - \Lambda(x) \rightarrow 0$. Thus, 
		$$\Lambda: (X,\tau_{\rho})\rightarrow ({\mathbb R},\tau_{|\cdot|})$$
		is sequentially continuous and therefore continuous, according to Theorem \ref{continuitysequential}.\\
		The implication $(iii)\Rightarrow (iv)$ is straightforward. Indeed, if $(iv)$ didn't hold, there would be a sequence $(x_j)$ with $\rho(x_j)<1$ and $|\Lambda(x_j)|>j$. Convexity would then imply $\rho\left(\frac{x_j}{j}\right)\rightarrow 0$ and $|\Lambda(x_j)|>1$, which contradicts $(iii)$. \\
		Assuming $(iv)$, if $\rho(x) = 1$, then $|\Lambda(2^{-1}x)|\leq \sup \limits_{\{x:\rho(x)<1\}} |\Lambda(x)|$. Hence, $|\Lambda(x)|\leq 2 \sup \limits_{\{x:\rho(x)<1\}} |\Lambda(x)|$ for any $x$ in the closed modular unit ball, i.e., $(v)$ holds. Conversely, assume $(v)$ holds and $r\geq 1$, then $\rho(x - x_0)\leq r$ implies, on account of convexity, $\rho\left(r^{-1}(x - x_0)\right)\leq 1$ and it follows that $|\Lambda(x)|\leq |\Lambda(x_0)| + r\sup\limits_{z \in \{x:\rho(x)\leq 1\}} |\Lambda(z)|$. On the other hand, if $0<r<1$ and $\rho(x - x_0) < r$, then $|\Lambda(x)|\leq |\Lambda(x_0)| + \sup\limits_{z \in \{x:\rho(x)\leq 1\}} |\Lambda(z)|$. Therefore, $(v)\Rightarrow (vi)$, and it is obvious that $(vi)\Rightarrow (vii)$.\\
		Any set $A\subset X$ with $\delta_{\rho}(A)=\sup\{\rho(a-b), a \in A, b \in A\}<\infty$ is contained in some modular ball, hence the implication $(vii)\Rightarrow (viii)$ is trivial. Finally, $(i)\Rightarrow (ix)$ follows directly from the fact that the modular topology of $X$ is contained in the norm topology of $X$.
	\end{proof}
	
	In the following example, we show that $X^{\rho}$ may be strictly contained in $X^{\ast}$, and hence $(ix)$ does not imply $(i)$. In particular, the boundedness of a linear functional on the modular unit ball does not imply modular continuity.
	
	\begin{example}\label{modular-example}
		{\normalfont 
			A simple analysis shows that one cannot expect any inclusion relation between the modular topology $\tau_{\rho}$ and the weak topology. To that end, take nonzero $\Lambda \in (\ell^{(n)})^* = (\ell^{\infty})^*$, with $\Lambda$ vanishing on the space of sequences that are convergent to $0$, denoted by ${\mathcal C}_0$. Let $(a_j) \in \ell^{(n)}$ be chosen so that $\Lambda((a_j)) = 1$. Since $(a_j) \in \ell^{(n)}$, there exists $\lambda \in (0,1)$ such that
			$$\rho((\lambda a_n)) = \textstyle \sum\limits_{n=1}^\infty \lambda^n |a_n|^n < \infty.$$
			Consider then the sequences $(c_n) = (\lambda a_j) {\mathbbm 1}_{j \geq n}$ and $(d_n) = (\lambda a_j){\mathbbm 1}_{j \leq n}$ in $\ell^{(n)}$.\\
			Since $\Lambda$ vanishes on ${\mathcal C}_0 = \ell^{\infty}_a$, one has
			$$\Lambda(0,0,...,0, \lambda a_n, \lambda a_{n+1},...) = \Lambda (\lambda a_1, \lambda a_2,...)\neq 0.$$
			Clearly,
			$$\lim\limits_{n \to \infty} \Lambda(c_n) = \lim\limits_{n \to \infty} \Lambda((\lambda a_j ){\mathbbm 1}_{j \geq n}) = \Lambda((\lambda a_j)) \neq 0,$$
			whereas 
			$$\lim\limits_{n \to \infty} \rho[(\lambda a_j){\mathbbm 1}_{j \geq n}] = 0.$$ 
			Modular convergence, then, does not imply weak convergence, and the weak topology of $\ell^{\infty}$ is not contained in the modular topology $\tau_{\rho}$. In particular, $\Lambda\in \left(\ell^{(n)}\right)^{\ast}\setminus \left(\ell^{(n)}\right)^{\rho}.$
		}
	\end{example}
	
	This example demonstrates that the weak topology and the modular topology are not directly comparable, as weak convergence may not align with modular convergence. In particular, it shows that even though a sequence may converge modularly, it may not converge weakly, indicating the absence of an inclusion relation between these topologies.
	\section{Applications and final remarks}\label{applications}
	The main driving force behind this investigation are the difficulties arising in the problem of minimization of the variable exponent $p(x)$-Dirichlet energy integral in the case when $p(x)$ finite a.e and unbounded. Specifically, let $\Omega\subset{\mathbb R}^n$ be a bounded, smooth domain and let $\int\limits_{\Omega}p^{-1}|\nabla \varphi|^pdx<\infty,$  with $n\leq\alpha <p(x)<\infty$ on $\Omega$. The classical minimization problem reads as follows: Minimize
	\begin{equation}\label{minimizationproblem}
		H \ni u\rightarrow 	F(u)=\int\limits_{\Omega}|\nabla (u-\varphi)(x)|^{p(x)}dx
	\end{equation}
	over a subspace $H$ of $W^{1,p(\cdot)}(\Omega)$ (the usual Sobolev space,\cite{DHHR,KR}) consisting of functions vanishing at the boundary. We refer the reader to \cite{DHHR,H, HHKV} and the references therein for the rather delicate discussion of the notion of zero boundary values of functions in $W^{1,p(\cdot)}(\Omega)$. \\
	If it exists, a solution $u$ to the minimization problem (\ref{minimizationproblem}) gives rise to a weak solution $w=\varphi-u$ of the Dirichlet problem
	\begin{equation}\label{DP}
		\begin{cases}
			\Delta_{p}w = \text{div}\left( |\nabla w|^{p-2} \nabla w \right) = 0 & \text{in} \ \Omega, \\
			w|_{\partial \Omega} = \varphi
		\end{cases}
	\end{equation}
	more precisely, the solution of the minimization problem yields a function $w$ such that 
	\begin{equation*}
		\int\limits_{\Omega}\left|\nabla w\right|^{p-2}\nabla w\nabla \phi \, dx=0\,\,\,\text{for any}\,\,\phi \in C^{\infty}_0(\Omega)
	\end{equation*}
	and $w-\varphi\in H$. The indispensability of modular topologies is revealed in the process of the natural minimization approach. It can be shown (though highly non-obvious) that any minimizing sequence $(u_j)$ of the Dirichlet energy integral is modularly convergent and that its modular limit is, in fact, the unique minimizer of the energy integral. Since the exponent $p(x)$ is unbounded in $\Omega$, theorem \ref{r-continuous-p^+L} implies that modular convergence is strictly weaker than norm convergence. This is the reason that the natural choice for $H$ is the modular closure of the test functions $C^{\infty}_0(\Omega)$ in $W^{1,p(\cdot)}(\Omega)$, which will be denoted by $V^{1,p(\cdot)}_0(\Omega)$. In $V^{1,p(\cdot)}_0(\Omega)$, thus, there is a unique minimizer of the Dirichlet integral and the Dirichlet problem (\ref{DP}) has a unique solution. We refer the interested reader to \cite{OAP, AOJA} for the detailed proofs and further discussions of the ideas barely sketched in this section.

\end{document}